\documentclass[a4paper,12pt,twoside]{article}
\usepackage{amsmath,amssymb,ifthen}
\usepackage[amsmath,thmmarks]{ntheorem}
\usepackage{paper}
\usepackage[all]{xy}
\usepackage{mathrsfs}
\SelectTips{cm}{11}
\newcommand{\rig}{\mathrm{rig}}
\newcommand{\red}{\mathrm{red}}
\newcommand{\M}{\breve{\mathscr{M}}}
\newcommand{\Mbar}{\bar{\mathscr{M}}}
\newcommand{\X}{\mathbb{X}}
\DeclareMathOperator{\Ind}{Ind}

\DeclareMathOperator{\Spf}{Spf}
\DeclareMathOperator{\Spa}{Spa}
\DeclareMathOperator{\GSp}{GSp}
\DeclareMathOperator{\rank}{rank}
\DeclareMathOperator{\supp}{supp}

\DeclareMathOperator{\Stab}{Stab}
\DeclareMathOperator{\Sh}{Sh}
\DeclareMathOperator{\Nilp}{\mathbf{Nilp}}

\DeclareMathOperator{\cInd}{c-Ind}

\title{Cuspidal representations in the $\ell$-adic cohomology of the Rapoport-Zink space for $\GSp(4)$}
\author{Tetsushi Ito and Yoichi Mieda}
\def\headertitle{$\ell$-adic cohomology of the Rapoport-Zink space for $\GSp(4)$}

\begin{document}

\maketitle

\begin{firstfootnote}
\end{firstfootnote}

\begin{abstract}
 In this paper, we study the $\ell$-adic cohomology of the Rapoport-Zink tower for GSp(4).
 We prove that the smooth representation of $\GSp_4(\Q_p)$ obtained as the $i$th compactly supported $\ell$-adic
 cohomology of the Rapoport-Zink tower has no quasi-cuspidal subquotient unless $i=2,3,4$.
 Our proof is purely local and does not require global automorphic methods.
\end{abstract}

\section{Introduction}
In \cite{MR1393439}, M.~Rapoport and Th.~Zink introduced certain moduli spaces of quasi-isogenies of
$p$-divisible groups with additional structures called the \textit{Rapoport-Zink spaces}.
They constructed systems of rigid analytic coverings of them which we call the \textit{Rapoport-Zink towers},
and established the $p$-adic uniformization theory of Shimura varieties generalizing classical
\v{C}erednik-Drinfeld uniformization. These spaces uniformize the rigid spaces associated with
the formal completion of certain Shimura varieties along Newton strata.

Using the $\ell$-adic cohomology of the Rapoport-Zink tower, we can construct a representation of
the product $G(\Q_p)\times J(\Q_p)\times W(\overline{\Q}_p/\Q_p)$, where $G$ is the reductive group
over $\Q_p$ corresponding to the Shimura datum, $J$ is an inner form of it, and $W(\overline{\Q}_p/\Q_p)$
is the Weil group of the $p$-adic field $\Q_p$. It is widely believed that this realizes
the local Langlands and Jacquet-Langlands correspondences (\cf \cite{MR1403942}).
Classical examples of the Rapoport-Zink spaces are the Lubin-Tate space and the Drinfeld upper half space;
these spaces were extensively studied by many people and many important results were obtained
(\cf \cite{MR0422290}, \cite{MR1044827}, \cite{MR1464867}, \cite{MR1876802}, \cite{MR2308851}, \cite{MR2511742}
and references therein). However, very little was known about the $\ell$-adic cohomology of other Rapoport-Zink
spaces.

The aim of this paper is to study cuspidal representations in the $\ell$-adic cohomology of the 
Rapoport-Zink tower for $\GSp_4(\Q_p)$.
Let us denote the Rapoport-Zink space for $\GSp_4(\Q_p)$ by $\M$.
It is a special formal scheme over $\Z_{p^\infty}=W(\overline{\F}_p)$ in the sense of Berkovich \cite{MR1395723}.
Let $\M^\rig$ be the Raynaud generic fiber of $\M$, that is, the generic fiber of
the adic space $t(\M)$ associated with $\M$. Using level structures at $p$, we can construct the
Rapoport-Zink tower
\[
 \cdots\longrightarrow \M^\rig_{m+1}\longrightarrow \M^\rig_m\longrightarrow \cdots\longrightarrow \M^\rig_2
 \longrightarrow \M^\rig_1\longrightarrow \M^\rig_0=\M^\rig,
\]
where $\M^\rig_m\longrightarrow \M^\rig$ is an \'etale Galois covering of rigid spaces
with Galois group $\GSp_4(\Z/p^m\Z)$.
We take the compactly supported $\ell$-adic cohomology (in the sense of \cite{MR1626021}) and take
the inductive limit of them. Then, on
\[
 H^i_{\mathrm{RZ}}:=\varinjlim_{m}H^i_c(\M^\rig_m\otimes_{\Q_{p^\infty}}\overline{\Q}_{p^\infty},\Q_\ell)
\]
(here $\Q_{p^\infty}=\Frac \Z_{p^\infty}$), we have an action of a product
\[
 \GSp_4(\Q_p)\times J(\Q_p)\times W(\overline{\Q}_p/\Q_p),
\]
where $J$ is an inner form of $\GSp_4$.

The main theorem of this paper is as follows:

\begin{thm}[Theorem \ref{thm:non-cusp-statement}]\label{thm:main-thm-introduction}
 The $\GSp_4(\Q_p)$-representation $H^i_{\mathrm{RZ}}\otimes_{\Q_\ell}\overline{\Q}_\ell$ has
 no quasi-cuspidal subquotient unless $i=2,3,4$.
\end{thm}
For the definition of quasi-cuspidal representations, see \cite[1.20]{MR771671}.
Note that since $\M^\rig_m$ is 3-dimensional for every $m\ge 0$, $H^i_{\mathrm{RZ}}=0$ unless $0\le i\le 6$.

Our proof of this theorem is purely local. We do not use global automorphic methods.
The main strategy of the proof is similar to that of \cite{non-cusp}, in which
the analogous result for the Lubin-Tate tower is given;
we construct the formal model $\M_m$ of $\M_m^\rig$ by using Drinfeld level structures and consider the geometry of
its special fiber. However, our situation is much more difficult
than the case of the Lubin-Tate tower. 
In the Lubin-Tate case, the tower consists of affine formal schemes $\{\Spf A_m\}_{m\ge 0}$,
and we can associate it with the tower of affine schemes $\{\Spec A_m\}_{m\ge 0}$. 
In \cite{non-cusp}, the second author defined the stratification on the special fiber of $\Spec A_m$
by using the kernel of the universal Drinfeld level structure,
and considered the local cohomology of the nearby cycle complex $R\psi\Lambda$ along the strata.
On the other hand, our tower $\{\M_m\}_{m\ge 0}$ does not consist of affine formal schemes
and there is no canonical way to associate it with a tower of schemes.
To overcome this problem, we take a sheaf-theoretic approach. For each direct summand $I$ of $(\Z/p^m\Z)^4$, 
we will define the complex of sheaves $\mathcal{F}_{m,I}$ on $(\M_m)_\red$ so that
the cohomology $H^i((\M_m)_\red,\mathcal{F}_{m,I})$ substitutes for the local cohomology of $R\psi\Lambda$ along
the strata defined by $I$ in the Lubin-Tate case. For the definition of $\mathcal{F}_{m,I}$,
we use the $p$-adic uniformization theorem by Rapoport and Zink.

There is another difficulty;
since a connected component of $\M$ is not quasi-compact,
the representation $H^i_{\mathrm{RZ}}$ of $\GSp_4(\Q_p)$ is far from admissible.
Therefore it is important to consider the action of $J(\Q_p)$ on $H^i_{\mathrm{RZ}}$, though
it does not appear in our main theorem. However, the cohomology $H^i((\M_m)_\red,\mathcal{F}_{m,I})$
has no apparent action of $J(\Q_p)$, since $J(\Q_p)$ does not act on the Shimura variety uniformized by $\M$.
We use the variants of formal nearby cycle 
introduced by the second author in \cite{formalnearby} to endow it with an action of $J(\Q_p)$.
Furthermore, to ensure the smoothness of this action, we use a property of finitely generated pro-$p$
groups (Section \ref{sec:pro-p}).
In fact, extensive use of the formalism developed in \cite{formalnearby} make us possible to
work mainly on the Rapoport-Zink tower itself and avoid the theory of $p$-adic uniformization
except for proving that $\M_m$ is locally algebraizable. However, for the reader's convenience,
we decided to make this article as independent of \cite{formalnearby} as possible.

The authors expect that the converse of Theorem \ref{thm:main-thm-introduction} also holds.
Namely, we expect that $H^i_{\mathrm{RZ}}\otimes_{\Q_\ell}\overline{\Q}_\ell$ has a quasi-cuspidal subquotient
if $i=2,3,4$. We hope to investigate it in a future work.

The outline of this paper is as follows. In Section \ref{sec:pro-p}, we prepare a criterion for 
the smoothness of representations over $\Q_\ell$. It is elementary but very powerful for our purpose.
In Section \ref{sec:Rapoport-Zink-space}, we give some basic definitions concerning with
the Rapoport-Zink space for $\GSp(4)$ and state the main theorem. 
Section \ref{sec:p-adic-unif} is devoted to introduce certain Shimura varieties related to our
Rapoport-Zink tower and recall the theory of $p$-adic uniformization.
The proof of the main theorem is accomplished in Section \ref{section:proof}.
The final Section \ref{sec:coh-corr} is an appendix on cohomological correspondences.
The results in the section are used to define actions of $\GSp_4(\Q_p)$ on various cohomology groups.

\bigbreak

\noindent{\bfseries Acknowledgment}\quad
The second author would like to thank Noriyuki Abe and Naoki Imai for the stimulating discussions.

\bigbreak

\noindent{\bfseries Notation}\quad
Let $p$ be a prime number and take another prime $\ell$ with $\ell\neq p$.
We denote the completion of the maximal unramified extension of $\Z_p$ by $\Z_{p^\infty}$ and its fraction field
by $\Q_{p^\infty}$.
Let $\Nilp=\Nilp_{\Z_{p^\infty}}$ be the category of $\Z_{p^\infty}$-schemes on which $p$ is locally nilpotent.
For an object $S$ of $\Nilp$, we put $\overline{S}=S\otimes_{\Z_{p^\infty}}\overline{\F}_p$.

In this paper, we use the theory of adic spaces (\cite{MR1306024}, \cite{MR1734903})
as a framework of rigid geometry.
A rigid space over $\Q_{p^\infty}$ is understood as an adic space locally of finite type over
$\Spa (\Q_{p^\infty},\Z_{p^\infty})$.

Every sheaf and cohomology are considered in the \'etale topology.
Every smooth representation is considered over $\Q_\ell$ or $\overline{\Q}_\ell$.
For a $\Q_\ell$-vector space $V$, we put $V_{\overline{\Q}_\ell}=V\otimes_{\Q_\ell}\overline{\Q}_\ell$.

\section{Preliminaries: smoothness of representations of profinite groups}\label{sec:pro-p}
Let $\mathbf{G}$ be a linear algebraic group over a $p$-adic field $F$.
In this section, we give a convenient criterion for the smoothness of a $\mathbf{G}(F)$-representation
over $\Q_\ell$.
The following theorem is essential:

\begin{thm}\label{thm:pro-p-main}
 Let $K$ be a closed subgroup of $\GL_n(\Z_p)$ and
 $(\pi,V)$ a finite-dimensional representation over $\Q_\ell$ of $K$ as an abstract group.
 Assume that there exists a $K$-stable $\Z_\ell$-lattice $\Lambda$ of $V$.
 Then this representation is automatically smooth.
\end{thm}

In order to prove this theorem, we require several facts on pro-$p$ groups.
Put $K_1=K\cap (1+p\mathrm{M}_n(\Z_p))$, which is a pro-$p$ open subgroup of $K$.

\begin{lem}
 The pro-$p$ group $K_1$ is (topologically) finitely generated.
\end{lem}

\begin{prf}
 By \cite[\S 5.1]{MR1720368}, the profinite group $\GL_n(\Z_p)$ has finite rank.
 In particular, $K_1$, a closed subgroup of $\GL_n(\Z_p)$, has finite topological generators.
\end{prf}

\begin{lem}\label{lem:strongly-complete}
 Every subgroup of finite index of $K_1$ is open.
\end{lem}

\begin{prf}
 In fact, this is true for every finitely generated pro-$p$ group;
 this is due to Serre \cite[4.2, Exercices 6)]{MR1324577}.
 See also \cite[Theorem 1.17]{MR1720368}, which gives a complete proof.
\end{prf}

\begin{rem}
 More generally, every subgroup of finite index of a finitely generated profinite group
 is open (\cite{MR2016979}, \cite{MR2276769}, \cite{MR2276770}). It is a very deep theorem.
\end{rem}

\begin{lem}\label{lem:hom-trivial}
 Let $G$ be a pro-$\ell$ group.
 Then every homomorphism $f\colon K_1\longrightarrow G$ is trivial.
\end{lem}

\begin{prf}
 Let $H$ be an open normal subgroup of $G$ and denote
 the composite $K_1\yrightarrow{f} G\longrightarrow G/H$ by $f_H$.
 By Lemma \ref{lem:strongly-complete}, $\Ker f_H$ is an open normal subgroup of $K_1$.
 Thus $K_1/\Ker f_H$ is a finite $p$-group.
 On the other hand, $G/H$ is a finite $\ell$-group.
 Since we have an injection $K_1/\Ker f_H\hooklongrightarrow G/H$, we have
 $K_1/\Ker f_H=1$, in other words, $f_H=1$.
 Therefore the composite $K_1\yrightarrow{f} G\yrightarrow{\cong} \varprojlim_{H} G/H$ is trivial.
 Hence we have $f=1$, as desired.
\end{prf}

\begin{prf}[of Theorem \ref{thm:pro-p-main}]
 Since $K_1$ is an open subgroup of $K$, we may replace $K$ by $K_1$.
 Take a $K_1$-stable $\Z_\ell$-lattice $\Lambda$ of $V$.
 Then, $\Lambda/\ell\Lambda$ is a finite abelian group.
 Therefore, by Lemma \ref{lem:strongly-complete}, there exists an open subgroup $U$ of $K_1$
 which acts trivially on $\Lambda/\ell\Lambda$. In other words,
 the homomorphism $\pi\colon K_1\longrightarrow \GL(\Lambda)\subset \GL(V)$ maps
 $U$ into the subgroup $1+\ell\End(\Lambda)$.
 Since $U$ is a closed subgroup of $1+p\mathrm{M}_n(\Z_p)$ and $1+\ell\End(\Lambda)$ is a pro-$\ell$ group,
 by Lemma \ref{lem:hom-trivial},
 the homomorphism $\pi\vert_U\colon U\longrightarrow 1+\ell\End(\Lambda)$ is trivial.
 Namely, $\pi\vert_U$ is a trivial representation.
\end{prf}

\begin{lem}\label{lem:alg-gp}
 Let $F$ be a $p$-adic field and $\mathbf{G}$ a linear algebraic group over $F$.
 Then every compact subgroup $K$ of $\mathbf{G}(F)$ can be realized as a closed subgroup of $\GL_n(\Z_p)$
 for some $n$.
\end{lem}

\begin{prf}
 Take an embedding $\mathbf{G}\hooklongrightarrow \GL_m$ defined over $F$.
 Since $\mathbf{G}(F)$ is a closed subgroup of $\GL_m(F)$, $K$ is also a closed subgroup of $\GL_m(F)$.
 Therefore we have a faithful continuous action of $K$ on $F^m$.
 By taking a $\Q_p$-basis of $F$, we have a faithful continuous action of $K$ on $\Q_p^n$
 for some $n$. Since $K$ is compact, it is well-known that there is a $K$-stable $\Z_p$-lattice
 in $\Q_p^n$. Hence we have a continuous injection $K\hooklongrightarrow \GL_n(\Z_p)$.
 Since $K$ is compact, it is isomorphic to a closed subgroup of $\GL_n(\Z_p)$.
\end{prf}

\begin{cor}\label{cor:smoothness-criterion}
 Let $F$ and $\mathbf{G}$ be as in the previous proposition.
 Let $I$ be a filtered ordered set and $\{K_i\}_{i\in I}$ be a system of compact open subgroups
 of $\mathbf{G}(F)$ indexed by $I$.

 Let $(\pi,V)$ be a (not necessarily finite-dimensional) $\Q_\ell$-representation of $\mathbf{G}(F)$ as an
 abstract group. Assume that there exists an inductive system $\{V_i\}_{i\in I}$ of
 finite-dimensional $\Q_\ell$-vector spaces satisfying the following:
 \begin{itemize}
  \item For every $i\in I$, $V_i$ is endowed with an action of $K_i$ as an abstract group.
  \item For every $i\in I$, $V_i$ has a $K_i$-stable $\Z_\ell$-lattice.
  \item There exists an isomorphism $\varinjlim_{i\in I}V_i\yrightarrow{\cong} V$ as 
	$\Q_\ell$-vector spaces such that the composite
	$V_i\longrightarrow \varinjlim_{i\in I}V_i\yrightarrow{\cong} V$ is $K_i$-equivariant
	for every $i\in I$.
 \end{itemize}
 Then $(\pi,V)$ is a smooth representation of $\mathbf{G}(F)$.
\end{cor}

\begin{prf}
 Let us take $x\in V$ and show that $\Stab_{\mathbf{G}(F)}(x)$, the stabilizer of $x$ in $\mathbf{G}(F)$,
 is open.
 There exists an element $i\in I$ such that $x$ lies in the image of
 $V_i\longrightarrow V$. Take $y\in V_i$ which is mapped to $x$.
 By Theorem \ref{thm:pro-p-main} and Lemma \ref{lem:alg-gp}, $V_i$ is a smooth representation
 of $K_i$. Therefore $\Stab_{K_i}(y)$ is open in $K_i$, hence is open in $\mathbf{G}(F)$.
 Since $V_i\longrightarrow V$ is
 $K_i$-equivariant, we have $\Stab_{K_i}(y)\subset \Stab_{K_i}(x)\subset \Stab_{\mathbf{G}(F)}(x)$.
 Thus $\Stab_{\mathbf{G}(F)}(x)$ is open in $\mathbf{G}(F)$, as desired.
\end{prf}

\begin{rem}
 Although we need the corollary above only for the case $F=\Q_p$, we proved it for a general $p$-adic field $F$
 for the completeness.
\end{rem}

\section{Rapoport-Zink space for $\GSp(4)$}\label{sec:Rapoport-Zink-space}
\subsection{The Rapoport-Zink space for $\GSp(4)$ and its rigid analytic coverings}
In this subsection, we recall basic definitions concerning with Rapoport-Zink spaces.
General definitions are given in \cite{MR1393439}, but here we restrict them to our special case.

Let $\X$ be a 2-dimensional isoclinic $p$-divisible group over $\overline{\F}_p$ with slope $1/2$,
and $\lambda_0\colon \X\yrightarrow{\cong} \X^\vee$ a (principal) polarization of $\X$, namely,
an isomorphism satisfying $\lambda_0^\vee=-\lambda_0$.
Consider the contravariant functor $\M\colon \Nilp\longrightarrow \mathbf{Set}$ that associates $S$
with the set of isomorphism classes of pairs $(X,\rho)$ consisting of
\begin{itemize}
 \item a 2-dimensional $p$-divisible group $X$ over $S$,
 \item and a quasi-isogeny (\cf \cite[Definition 2.8]{MR1393439})
       $\rho\colon \X\otimes_{\overline{\F}_p}\overline{S}\longrightarrow X\otimes_S\overline{S}$,
\end{itemize}
such that there exists an isomorphism $\lambda\colon X\longrightarrow X^\vee$ which makes the following
diagram commutative up to multiplication by $\Q_p^\times$:
	\[
	 \xymatrix{%
	 \mathbb{X}\otimes_{\overline{\F}_p}\overline{S}\ar[r]^-{\rho}\ar[d]^-{\lambda_0\otimes\id}& X\otimes_S\overline{S}\ar[d]^-{\lambda\otimes \id}\\
	 \mathbb{X}^\vee\otimes_{\overline{\F}_p}\overline{S}& X^\vee\otimes_S\overline{S}\lefteqn{.}\ar[l]_-{\rho^\vee}
	}
	\]
Note that such $\lambda$ is uniquely determined by $(X,\rho)$ up to multiplication by $\Z_p^\times$
and gives a polarization of $X$.
It is proved by Rapoport-Zink that $\M$ is represented by a special formal scheme (\cf \cite{MR1395723})
over $\Spf \Z_{p^\infty}$.
Moreover, $\M$ is separated over $\Spf \Z_{p^\infty}$
\cite[Lemme 2.3.23]{MR2074714}. However, $\M$ is neither quasi-compact nor $p$-adic.
We put $\Mbar=\M_{\mathrm{red}}$, which is a scheme locally of finite type and separated over $\overline{\F}_p$.
It is known that $\Mbar$ is $1$-dimensional (for example, see \cite{MR2466182}) and
every irreducible component of $\Mbar$ is projective over $\overline{\F}_p$ \cite[Proposition 2.32]{MR1393439}.
In particular, $\Mbar$ has a locally finite quasi-compact open covering.

Let $D(\X)_\Q=(N,\Phi)$ be the rational Dieudonn\'e module of $\X$, which is a 4-dimensional isocrystal
over $\Q_{p^\infty}$. The fixed polarization $\lambda_0$ gives the alternating pairing
$\langle\ ,\ \rangle_{\lambda_0}\colon N\times N\longrightarrow \Q_{p^\infty}(1)$.
We define the algebraic group $J$ over $\Q_p$ as follows: for a $\Q_p$-algebra $R$, the group $J(R)$
consists of elements $g\in \GL(R\otimes_{\Q_p}N)$ such that
\begin{itemize}
 \item $g$ commutes with $\Phi$,
 \item and $g$ preserves the pairing $\langle\ ,\ \rangle_{\lambda_0}$ up to scalar multiplication,
       i.e., there exists $c(g)\in R^\times$ such that 
       $\langle gx,gy\rangle_{\lambda_0}=c(g)\langle x,y\rangle_{\lambda_0}$ for every $x,y\in R\otimes_{\Q_p}N$.
\end{itemize}
It is an inner form of $\GSp(4)$, since $D(\X)_\Q$ is the isocrystal associated with
a basic Frobenius conjugacy class of $\GSp(4)$.

In the sequel, we also denote $J(\Q_p)$ by $J$.
Every element $g\in J$ naturally induces a quasi-isogeny $g\colon \X\longrightarrow \X$ and
the following diagram is commutative up to $\Q_p^\times$-multiplication:
\[
 \xymatrix{%
 \mathbb{X}\ar[r]^-{g}\ar[d]^-{\lambda_0}&\mathbb{X}\ar[d]^-{\lambda_0}\\
 \mathbb{X}^\vee&\mathbb{X}^\vee\lefteqn{.}\ar[l]_-{g^{\vee}}
 }
\]
Therefore, we can define the left action of $J$ on $\M$ by $g\colon \M(S)\longrightarrow \M(S)$;
$(X,\rho)\longmapsto (X,\rho\circ g^{-1})$.

We denote the Raynaud generic fiber of $\M$ by $\M^\rig$.
It is defined as $t(\M)\setminus V(p)$, where $t(\M)$ is the adic space associated with $\M$
(\cf \cite[Proposition 4.1]{MR1306024}). 
As $\M$ is separated and special over $\Z_{p^\infty}$, 
$\M^\rig$ is separated and locally of finite type over $\Spa (\Q_{p^\infty},\Z_{p^\infty})$.
Since $\M$ has a locally finite quasi-compact open covering, $\M^\rig$ is taut
by \cite[Lemma 4.14]{formalnearby}.
Moreover, by using the period morphism \cite[Chapter 5]{MR1393439}, we can see that $\M^\rig$ is
3-dimensional and smooth over $\Spa (\Q_{p^\infty},\Z_{p^\infty})$ (\cf \cite[Proposition 5.17]{MR1393439}).

Next we will consider level structures. 
Let $\widetilde{X}$ be the universal $p$-divisible group over $\M$
and $\widetilde{X}^\rig$ be the associated $p$-divisible group over $\M^\rig$.
Note that $\widetilde{X}^\rig$ is an \'etale $p$-divisible group.
Let us fix a polarization $\widetilde{\lambda}\colon \widetilde{X}\longrightarrow \widetilde{X}^\vee$ which is
compatible with $\lambda_0$, i.e., satisfies the condition in the definition of $\M$.
Let $\mathsf{S}$ be a connected rigid space over $\Q_{p^\infty}$ (i.e., a connected adic space locally of
finite type over $\Spa (\Q_{p^\infty},\Z_{p^\infty})$), $\mathsf{S}\longrightarrow \M^\rig$ a morphism
over $\Q_{p^\infty}$ and $\widetilde{X}^\rig_{\mathsf{S}}$ the pull-back of $\widetilde{X}^\rig$.
Fix a geometric point $\overline{x}$ of $\mathsf{S}$ and 
an isomorphism $T_p(\mu_{p^\infty,\mathsf{S}})_{\overline{x}}=\Z_p(1)\cong \Z_p$.
Then $\widetilde{\lambda}$ induces an alternating bilinear form $\psi_{\widetilde{\lambda}}$
on the $\pi_1(\mathsf{S},\overline{x})$-module $(T_p\widetilde{X}^\rig)_{\overline{x}}$;
\[
 \psi_{\widetilde{\lambda}}\colon (T_p\widetilde{X}^\rig)_{\overline{x}}\times (T_p\widetilde{X}^\rig)_{\overline{x}}\longrightarrow T_p(\mu_{p^\infty,\mathsf{S}})_{\overline{x}}\cong \Z_p.
\]
Fix a free $\Z_p$-module $L$ of rank 4 and a perfect alternating bilinear form
$\psi_0\colon L\times L\longrightarrow \Z_p$. Put $K_0=\GSp(L,\psi_0)$, $V=L\otimes_{\Z_p}\Q_p$ and
$G=\GSp(V,\psi_0)$.
Let $T(\mathsf{S},\overline{x})$ be the set consisting of isomorphisms $\eta\colon L\yrightarrow{\cong} (T_p\widetilde{X}^\rig)_{\overline{x}}$ which map $\psi_0$ to $\Z_p^\times$-multiples of $\psi_{\widetilde{\lambda}}$.
It is independent of the choice of $\widetilde{\lambda}$ and
$T_p(\mu_{p^\infty,\mathsf{S}})_{\overline{x}}\cong \Z_p$, since they are unique up to
$\Z_p^\times$-multiplication. Obviously, the groups $K_0$ and $\pi_1(\mathsf{S},\overline{x})$ naturally act on
$T(\mathsf{S},\overline{x})$. 

For an open subgroup $K$ of $K_0$, a $K$-level structure of $\widetilde{X}^\rig_\mathsf{S}$
means an element of $(T(\mathsf{S},\overline{x})/K)^{\pi_1(\mathsf{S},\overline{x})}$.
Note that, if we change a geometric point $\overline{x}$ to $\overline{x}'$, the sets
$(T(\mathsf{S},\overline{x})/K)^{\pi_1(\mathsf{S},\overline{x})}$ and
$(T(\mathsf{S},\overline{x}')/K)^{\pi_1(\mathsf{S},\overline{x}')}$ are naturally isomorphic.
Thus the notion of $K$-level structures is independent of the choice of $\overline{x}$.
The functor that associates $\mathsf{S}$ with the set of 
$K$-level structures of $\widetilde{X}^\rig_\mathsf{S}$ is represented by a finite Galois 
\'etale covering $\M^\rig_K\longrightarrow \M^\rig$, whose Galois group is $K_0/K$.
Since $T(\mathsf{S},\overline{x})$ is a $K_0$-torsor, $\M^\rig_{K_0}$ coincides with $\M^\rig$.
If $K'$ is an open subgroup of $K$, we have a natural morphism 
$p_{KK'}\colon \M^\rig_{K'}\longrightarrow \M^\rig_K$. Therefore, we get the projective system of rigid spaces
$\{\M^\rig_K\}_K$ indexed by the filtered ordered set of open subgroups of $K_0$,
which is called the \textit{Rapoport-Zink tower}.
Obviously, the group $J$ acts on the projective system $\{\M^\rig_K\}_K$.

Let $g$ be an element of $G$ and $K$ an open subgroup of $K_0$ which is enough small so that
$g^{-1}Kg\subset K_0$. Then we have a natural morphism
$\M^\rig_K\longrightarrow \M^\rig_{g^{-1}Kg}$ over $\Q_{p^\infty}$.
If $g\in K_0$, then it is given by $\eta\longmapsto \eta\circ g$;
for other $g$, it is more complicated \cite[5.34]{MR1393439}.
In any case, we get a right action of $G$ on the pro-object $\sideset{\text{``}}{\text{''}}\varprojlim\M^\rig_K$.

\begin{defn}
 We put $H^i_{\mathrm{RZ}}=\varinjlim_{K}H^i_c(\M^\rig_K\otimes_{\Q_{p^\infty}}\overline{\Q}_{p^\infty},\Q_\ell)$.
\end{defn}

Here $H^i_c(\M^\rig_K\otimes_{\Q_{p^\infty}}\overline{\Q}_{p^\infty},\Q_\ell)$
 is the compactly supported $\ell$-adic
cohomology of $\M^\rig_K\otimes_{\Q_{p^\infty}}\overline{\Q}_{p^\infty}$ defined in \cite{MR1626021};
note that $\M^\rig_K$ is separated and taut.
By the constructions above, $G\times J$ acts on $H^i_{\mathrm{RZ}}$ on the left
(the action of $j\in J$ is given by $(j^{-1})^*$).
Obviously the action of $G$ on $H^i_{\mathrm{RZ}}$ is smooth. On the other hand,
it is known that the action of $J$ on $H^i_{\mathrm{RZ}}$ is also smooth.
This is due to Berkovich (see \cite[Corollaire 4.4.7]{MR2074714}); see also Remark \ref{rem:Berkovich-smoothness},
where we give another proof of the smoothness.
Hence we get the smooth representation $H^i_{\mathrm{RZ}}$ of $G\times J$.

Our main theorem is the following:

\begin{thm}[Non-cuspidality]\label{thm:non-cusp-statement}
 The smooth representation $H^i_{\mathrm{RZ},\overline{\Q}_\ell}$ of $G$ has no quasi-cuspidal subquotient unless $i=2,3,4$.
\end{thm}
For the definition of quasi-cuspidal representations, see \cite[1.20]{MR771671}.

Theorem \ref{thm:non-cusp-statement} is proved in Section \ref{section:proof}.

\subsection{An integral model $\M_m$ of $\M^\rig_{K_m}$}
For an integer $m\ge 1$, let $K_m$ be the kernel of $\GSp(L,\psi_0)\longrightarrow \GSp(L/p^mL,\psi_0)$.
It is an open subgroup of $K_0$. We can describe the definition of $K_m$-level structures more concretely.
As in the previous subsection, we fix a polarization $\widetilde{\lambda}$ of $\widetilde{X}^\rig$ which is
compatible with $\lambda_0$. It induces the alternating bilinear morphism between finite \'etale group schemes
$\psi_{\widetilde{\lambda}}\colon \widetilde{X}^\rig[p^m]\times \widetilde{X}^\rig[p^m]\longrightarrow \mu_{p^m}$.
Let $\mathsf{S}\longrightarrow \M^\rig$ be as in the previous subsection.
Then a $K_m$-level structure of $\widetilde{X}^\rig_\mathsf{S}$ naturally corresponds bijectively
to an isomorphism $\eta\colon L/p^mL\yrightarrow{\cong}\widetilde{X}_{\mathsf{S}}^\rig[p^m]$ between finite
\'etale group schemes such that there exists an isomorphism $\Z/p^m\Z\yrightarrow{\cong}\mu_{p^m,\mathsf{S}}$
which makes the following diagram commutative:
\[
 \xymatrix{%
 L/p^mL\times L/p^mL\ar[r]^-{\psi_0}\ar[d]_-{\eta\times\eta}^{\cong}& \Z/p^m\Z\ar[d]^-{\cong}\\
 \widetilde{X}^\rig_{\mathsf{S}}[p^m]\times \widetilde{X}^\rig_{\mathsf{S}}[p^m]\ar[r]^-{\psi_{\widetilde{\lambda}}}& \mu_{p^m,\mathsf{S}}\lefteqn{.}
 }
\]
For simplicity, we write $\M_m^\rig$ for $\M_{K_m}^\rig$ and $p_{mn}$ for $p_{K_mK_n}$.
In this subsection, we construct a formal model $\M_m$ of $\M_m^\rig$ by following \cite[\S 6]{MR2169874}.
Let $\mathcal{S}$ be a formal scheme of finite type over $\M^\rig$ and denote by $\widetilde{X}_{\mathcal{S}}$ the pull-back
of $\widetilde{X}$ to $\mathcal{S}$. A Drinfeld $m$-level structure of $\widetilde{X}_{\mathcal{S}}$
is a morphism $\eta\colon L/p^mL\longrightarrow \widetilde{X}_{\mathcal{S}}[p^m]$ satisfying
the following conditions:
\begin{itemize}
 \item the image of $\eta$ gives a full set of sections of $\widetilde{X}_{\mathcal{S}}[p^m]$,
 \item and there exists a morphism $\Z/p^m\Z\longrightarrow \mu_{p^m,\mathcal{S}}$ 
       which makes the following diagram commutative:
       \[
       \xymatrix{%
       L/p^mL\times L/p^mL\ar[r]^-{\psi_0}\ar[d]_-{\eta\times\eta}& \Z/p^m\Z\ar[d]\\
       \widetilde{X}_{\mathcal{S}}[p^m]\times \widetilde{X}_{\mathcal{S}}[p^m]\ar[r]^-{\psi_{\widetilde{\lambda}}}& \mu_{p^m,\mathcal{S}}\lefteqn{.}
       }
       \]
\end{itemize}
It is known that the functor that associates $\mathcal{S}$ with the set of 
Drinfeld $m$-level structures of $\widetilde{X}_\mathcal{S}$ is represented by the formal scheme $\M_m$
which is finite over $\M$ (\cf \cite[Proposition 15]{MR2169874}).
Note that, unlike the case of Lubin-Tate tower, $\M_m$ is not necessarily flat over $\M$.
It is easy to show that $\M_m$ gives a formal model of $\M^\rig_m$, namely,
the Raynaud generic fiber of $\M_m$ coincides with $\M^\rig_m$.
We denote $(\M_m)_{\mathrm{red}}$ by $\Mbar_m$, which is a 1-dimensional scheme over $\overline{\F}_p$.

There is a natural left action of $J$ on $\M_m$ which is compatible with
that on $\M^\rig_m$. On the other hand, the natural action $K_0$ on $L/p^mL$ induces a right action
of $K_0$ on $\M_m$, which is compatible with that on $\M_{K_m}^\rig$. 

We can also describe $\M_m$ as a functor from $\Nilp$ to $\mathbf{Set}$; for an object $S$ of $\Nilp$,
the set $\M_m(S)$ consists of isomorphism classes of triples $(X,\rho,\eta)$,
where $(X,\rho)\in \M_m(S)$ and $\eta\colon L/p^mL\longrightarrow X[p^m]$ is a Drinfeld $m$-level structure
of $X$. By this description, the action of $j\in J$ on $\M_m$ is given by 
$(X,\rho,\eta)\longmapsto (X,\rho\circ j^{-1},\eta)$. On the other hand, the action of $g\in K_0$ on $\M_m$
is given by $(X,\rho,\eta)\longmapsto (X,\rho,\eta\circ g)$.

By \cite[Lemma 7.2]{MR2074715}, $\{\M_m\}_{m\ge 0}$ forms a projective system of formal schemes
equipped with the commuting action of $J$ and $K_0$.

\subsection{Compactly supported cohomology of $\Mbar_m$}
For $m\ge 0$, we denote the set of quasi-compact open subsets of $\Mbar_m$ by $\mathcal{Q}_m$.
It has a natural filtered order by inclusion.

\begin{defn}
 For an object $\mathcal{F}$ of $D^b(\Mbar_m,\Z_\ell)$ or $D^b(\Mbar_m,\Q_\ell)$,
 we put 
 \[
 H^i_c(\Mbar_m,\mathcal{F})=\varinjlim_{U\in\mathcal{Q}_m}H^i_c(U,\mathcal{F}\vert_U).
 \]
 Assume that $\mathcal{F}$ has a $J$-equivariant structure, namely, for every $g\in J$ an isomorphism
 $\varphi_g\colon g^*\mathcal{F}\yrightarrow{\cong} \mathcal{F}$ is given such that
 $\varphi_{gg'}=\varphi_{g'}\circ g'^*\varphi_g$ for every $g,g'\in J$.
 Then $J$ naturally acts on $H^i_c(\Mbar_m,\mathcal{F})$ on the right.
 Therefore we get a left action of $J$ on $H^i_c(\Mbar_m,\mathcal{F})$ by taking the inverse
 $J\longrightarrow J$; $g\longmapsto g^{-1}$.
\end{defn}

\begin{thm}\label{thm:J-fg-smooth}
 Let $\mathcal{F}^\circ$ be an object of $D_c^b(\Mbar_m,\Z_\ell)$
 and $\mathcal{F}$ the object of $D_c^b(\Mbar_m,\Q_\ell)$ associated with $\mathcal{F}^\circ$.
 Assume that we are given a $J$-equivariant structure of $\mathcal{F}^\circ$ 
 (thus $\mathcal{F}$ also has a $J$-equivariant structure).
 Then $H^i_c(\Mbar_m,\mathcal{F})$ is a finitely generated smooth $J$-representation.
\end{thm}

\begin{prf}
 Let $U$ be an element of $\mathcal{Q}_m$. By \cite[Proposition 2.3.11]{MR2074714},
 there exists a compact open subgroup $K_U$ of $J$ which stabilizes $U$.
 Then $H^i_c(U,\mathcal{F}\vert_U)$ is a finite-dimensional $\Q_\ell$-vector space endowed with the action of $K_U$
 and has the $K_U$-stable $\Z_\ell$-lattice
 $\Imm(H^i_c(U,\mathcal{F}^\circ\vert_U)\longrightarrow H^i_c(U,\mathcal{F}\vert_U))$.
 Therefore $H^i_c(\Mbar_m,\mathcal{F})$ is a smooth $J$-representation
 by Corollary \ref{cor:smoothness-criterion}.

 To prove that $H^i_c(\Mbar_m,\mathcal{F})$ is finitely generated,
 we may assume $m=0$, for $H^i_c(\Mbar_m,\mathcal{F})=H^i_c(\Mbar_0,p_{0m*}\mathcal{F})$.
 In this case, we can use the similar method as in \cite[Proposition 4.4.13]{MR2074714}.
 Let us explain the argument briefly. By \cite[Th\'eor\`eme 2.4.13]{MR2074714},
 there exists $W\in\mathcal{Q}_0$ such that $\bigcup_{g\in J}gW=\Mbar_0$.
 We put $K=\{g\in J\mid gW=W\}$ and $\Omega=\{g\in J\mid gW\cap W\neq\varnothing\}$.
 As in the proof of \cite[Proposition 4.4.13]{MR2074714}, $K$ is a compact open subgroup of $J$ and
 $\Omega$ is a compact subset of $J$. 
 For $\alpha=([g_1],\ldots,[g_n])\in (J/K)^n$, we put $W_\alpha=g_1W\cap \cdots \cap g_nW$
 and $K_\alpha=\bigcap_{j=1}^ng_jKg_j^{-1}$.
 For an open covering $\{gW\}_{g\in J/K}$, we can associate the \v{C}ech spectral sequence
 \[
 E_1^{r,s}=\bigoplus_{\alpha\in (J/K)^{-r+1}}H^s_c(W_\alpha,\mathcal{F}\vert_{W_\alpha})\Longrightarrow
 H^{r+s}_c(\Mbar_0,\mathcal{F}).
 \]
 Consider the diagonal action of $J$ on $(J/K)^{-r+1}$.
 The coset 
 \[
 J\backslash\{\alpha\in (J/K)^{-r+1}\mid W_\alpha\neq \varnothing\}
 \]
 is finite; indeed,
 if $W_\alpha\neq \varnothing$ for $\alpha=([g_1],\ldots,[g_{-r+1}])\in (J/K)^{-r+1}$,
 then $g_1^{-1}\alpha \in \{1\}\times \Omega/K\times\cdots\times \Omega/K$,
 which is a finite set.

 Take a system of representatives $\alpha_1,\ldots,\alpha_n$ of the coset above.
 Then there is a natural isomorphism
 $\bigoplus_{\alpha\in J\alpha_j}H^s_c(W_\alpha,\mathcal{F}\vert_{W_\alpha})\cong \cInd^{J}_{K_{\alpha_j}}H^s_c(W_{\alpha_j},\mathcal{F}\vert_{W_{\alpha_j}})$.
 Hence $E_1^{r,s}\cong \bigoplus_{j=1}^n\cInd^{J}_{K_{\alpha_j}}H^s_c(W_{\alpha_j},\mathcal{F}\vert_{W_{\alpha_j}})$ is a finitely generated $J$-module,
 since the cohomology $H^s_c(W_{\alpha_j},\mathcal{F}\vert_{W_{\alpha_j}})$ is finite-dimensional
 for each $j$.
 By this and the fact that a finitely generated smooth $J$-module is noetherian
 \cite[Remarque 3.12]{MR771671}, we conclude that $H^i_c(\Mbar_0,\mathcal{F})$ is finitely generated.
\end{prf}

\begin{lem}\label{lem:fixed-K_m}
 Let $\mathcal{F}$ be an object of $D^b_c(\Mbar_m,\Q_\ell)$ with a $K_0/K_m$-equivariant structure.
 Let $n$ be an integer with $0\le n\le m$ and put $\mathcal{G}=(p_{nm*}\mathcal{F})^{K_n/K_m}$.
 Then we have $H^i_c(\Mbar_m,\mathcal{F})^{K_n/K_m}=H^i_c(\Mbar_n,\mathcal{G})$.
\end{lem}

\begin{prf}
 Since the cardinality of $K_n/K_m$ is prime to $\ell$, $(-)^{K_n/K_m}$ commutes with $H^i_c$.
 Therefore, we have
 \begin{align*}
  &H^i_c(\Mbar_m,\mathcal{F})^{K_n/K_m}
  =\varinjlim_{U\in \mathcal{Q}_m} H^i_c(U,\mathcal{F}\vert_U)^{K_n/K_m}
  =\varinjlim_{V\in \mathcal{Q}_n} H^i_c\bigl(p_{nm}^{-1}(V),\mathcal{F}\vert_{p_{nm}^{-1}(V)}\bigr)^{K_n/K_m}\\
  &\qquad=\varinjlim_{V\in \mathcal{Q}_n} H^i_c\bigl(V,p_{nm*}(\mathcal{F}\vert_{p_{nm}^{-1}(V)})\bigr)^{K_n/K_m}
  =\varinjlim_{V\in \mathcal{Q}_n} H^i_c\Bigl(V,\bigl(p_{nm*}(\mathcal{F}\vert_{p_{nm}^{-1}(V)})\bigr)^{K_n/K_m}\Bigr)\\
  &\qquad =\varinjlim_{V\in \mathcal{Q}_n} H^i_c(V,\mathcal{G}\vert_V)=H^i_c(\Mbar_n,\mathcal{G}).
 \end{align*}
\end{prf}

\begin{defn}
 A \textit{system of coefficients} over the tower $\{\Mbar_m\}_{m\ge 0}$ is the data 
 $\mathcal{F}=\{\mathcal{F}_m\}_{m\ge 0}$ where $\mathcal{F}_m$ is an object of $D^b_c(\Mbar_m,\Q_\ell)$
 with a $K_0/K_m$-equivariant structure
 such that $(p_{nm*}\mathcal{F})^{K_n/K_m}=\mathcal{F}_n$ for every integers $m$, $n$ with $0\le n\le m$.
 Then, by Lemma \ref{lem:fixed-K_m}, we have 
 $H^i_c(\Mbar_m,\mathcal{F}_m)^{K_n/K_m}=H^i_c(\Mbar_n,\mathcal{F}_n)$.
 We put $H^i_c(\Mbar_\infty,\mathcal{F})=\varinjlim_{m}H^i_c(\Mbar_m,\mathcal{F}_m)$.

 If each $\mathcal{F}_m$ is endowed with a $J$-equivariant structure which commutes with the given
 $K_0/K_m$-equivariant structure,
 and for every $0\le n\le m$ the $J$-equivariant structures on $\mathcal{F}_m$ and $\mathcal{F}_n$ are compatible under the identification 
 $(p_{nm*}\mathcal{F}_m)^{K_n/K_m}=\mathcal{F}_n$, then we say that we have
 a $J$-equivariant structure on $\mathcal{F}$. Such a structure naturally induces the action of $J$
 on $H^i_c(\Mbar_\infty,\mathcal{F})$.

 By replacing ``$D^b_c(\Mbar_m,\Q_\ell)$'' with ``$D^b_c(\Mbar_m,\Z_\ell)$'',
 we may also define a \textit{system of integral coefficients} $\mathcal{F}^\circ$ over $\{\Mbar_m\}_{m\ge 0}$,
 the cohomology $H^i_c(\Mbar_\infty,\mathcal{F}^\circ)$ and a $J$-equivariant structure on $\mathcal{F}^\circ$.
\end{defn}

\begin{cor}\label{cor:M_infty-coh}
 Let $\mathcal{F}^\circ$ be a system of integral coefficients over $\{\Mbar_m\}_{m\ge 0}$
 with a $J$-equivariant structure and
 $\mathcal{F}$ the system of coefficients associated with $\mathcal{F}^\circ$.
 Then $H^i_c(\Mbar_\infty,\mathcal{F})$ is a smooth $K_0\times J$-representation and 
 $H^i_c(\Mbar_\infty,\mathcal{F})^{K_m}$ is a finitely generated smooth $J$-representation
 for every integer $m\ge 0$.
\end{cor}

\begin{prf}
 The smoothness is clear from Theorem \ref{thm:J-fg-smooth} and the definition of 
 $H^i_c(\Mbar_\infty,\mathcal{F})$.
 Since $H^i_c(\Mbar_\infty,\mathcal{F})^{K_m}=H^i_c(\Mbar_m,\mathcal{F}_m)$, the second assertion also follows
 from Theorem \ref{thm:J-fg-smooth}.
\end{prf}

\section{Shimura variety and $p$-adic uniformization}\label{sec:p-adic-unif}
In this section, we introduce certain Shimura varieties (Siegel threefolds) related to our Rapoport-Zink tower.
Let us fix a $4$-dimensional $\Q$-vector space $V'$ and an alternating perfect pairing 
$\psi'\colon V'\times V'\longrightarrow \Q$. For an integer $m\ge 0$ and a compact open subgroup
$K^p\subset \GSp(V'_{\A^{\infty,p}})=\GSp(V'_{\A^{\infty,p}},\psi'_{\A^{\infty,p}})$,
 consider the functor
$\Sh_{m,K^p}$ from the category of locally noetherian $\Z_{p^\infty}$-schemes to the category of sets
that associates $S$ with the set of isomorphism classes of quadruples $(A,\lambda,\eta^p,\eta_p)$ where
\begin{itemize}
 \item $A$ is a projective abelian surface over $S$ up to prime-to-$p$ isogeny,
 \item $\lambda\colon A\longrightarrow A^\vee$ is a prime-to-$p$ polarization,
 \item $\eta^p$ is a $K^p$-level structure of $A$,
 \item and $\eta_p\colon L/p^mL\longrightarrow A[p^m]$ is a Drinfeld $m$-level structure
\end{itemize}
(for the detail, see \cite[\S 5]{MR1124982}).
Two quadruples $(A,\lambda,\eta^p,\eta_p)$ and $(A',\lambda',\eta'^p,\eta'_p)$ are said to be isomorphic
if there exists a prime-to-$p$ isogeny from $A$ to $A'$ which carries $\lambda$ to a $\Z_{(p)}^\times$-multiple
of $\lambda'$, $\eta^p$ to $\eta'^p$ and $\eta_p$ to $\eta'_p$. 
We put $\Sh_{K^p}=\Sh_{0,K^p}$.
It is known that if $K^p$ is sufficiently small,
$\Sh_{m,K^p}$ is represented by a quasi-projective scheme over $\Z_{p^\infty}$ with smooth generic fiber. 
In the sequel, we always assume that $K^p$ is enough small so that $\Sh_{m,K^p}$ is representable.
We denote the special fiber of $\Sh_{m,K^p}$ (resp.\ $\Sh_{K^p}$) by $\overline{\Sh}_{m,K^p}$
(resp.\ $\overline{\Sh}_{K^p}$).

For a compact open subgroup $K'^p$ contained in $K^p$ and an integer $m'\ge m$, we have the natural
morphism $\Sh_{m',K'^p}\longrightarrow \Sh_{m,K^p}$. This is a finite morphism and is moreover \'etale
if $m'=m$.

Next we recall the $p$-adic uniformization theorem, which gives a relation between 
$\M$ and $\Sh_{K^p}$.
Let us fix a polarized abelian surface $(A_0,\lambda_{A_0})$ over $\overline{\F}_p$ such that
$A_0[p^\infty]$ is an isoclinic $p$-divisible group with slope $1/2$.
Note that such $(A_0,\lambda_{A_0})$ exists; 
for example, we can take $(A_0,\lambda_{A_0})=(E^2,\lambda_E^2)$,
where $E$ is a supersingular elliptic curve over $\overline{\F}_p$ and
$\lambda_E$ is a polarization of $E$. By definition, the rational Dieudonn\'e module
$D(A_0[p^\infty])_{\Q}$ is isomorphic to $D(\mathbb{X})_{\Q}$. Thus, by the subsequent lemma,
there is an isomorphism of isocrystals $D(A_0[p^\infty])_{\Q}\cong D(\mathbb{X})_{\Q}$ which preserves
the natural polarizations.

\begin{lem}\label{lem:RR}
 We use the notation in \cite[\S 1]{MR1411570}. Let $d\ge 1$ be an integer.
 \begin{enumerate}
  \item Let $b$ be an element of $B(\GSp_{2d})$ and $b'$ the image of $b$ under the natural map 
	$B(\GSp_{2d})\longrightarrow B(\GL_{2d})$.
	Then $b$ is basic if and only if $b'$ is basic.
  \item The map $B(\GSp_{2d})_{\mathrm{basic}}\longrightarrow B(\GL_{2d})_{\mathrm{basic}}$
	induced from i) is an injection.
 \end{enumerate}
\end{lem}

\begin{prf}
 Note that the center of $\GSp_{2d}$ coincides with that of $\GL_{2d}$. Thus i) is clear, 
 since $b$ (resp.\ $b'$) is basic if and only if the slope morphism $\nu_b\colon \mathbb{D}\longrightarrow \GSp_{2d}$ (resp.\ $\nu_{b'}\colon \mathbb{D}\yrightarrow{\nu_b}\GSp_{2d}\hooklongrightarrow \GL_{2d}$)
factors through the center of $\GSp_{2d}$ (resp.\ $\GL_{2d}$).

 We prove ii). By \cite[Theorem 1.15]{MR1411570}, it suffices to show that the natural map
 $\pi_1(\GSp_{2d})\longrightarrow \pi_1(\GL_{2d})$ is injective. Take a maximal torus $T$ (resp.\ $T'$) of
 $\GSp_{2d}$ (resp.\ $\GL_{2d}$) such that $T\subset T'$. Then, since $\mathrm{Sp}_{2d}$
 (resp.\ $\mathrm{SL}_{2d}$) is simply connected, $\pi_1(\GSp_{2d})$ (resp.\ $\pi_1(\GL_{2d})$)
 can be identified with the quotient of $X_*(T)$ (resp.\ $X_*(T')$) 
 induced by $c\colon T\twoheadlongrightarrow \mathbb{G}_m$ (resp.\ $\det\colon T'\twoheadlongrightarrow \mathbb{G}_m$), 
 where $c$ denotes the similitude character of $\GSp_{2d}$. In particular, both $\pi_1(\GSp_{2d})$ and
 $\pi_1(\GL_{2d})$ are isomorphic to $\Z$.

 The commutative diagram
 \[
  \xymatrix{%
 \GSp_{2d}\ar@{->>}[r]^-{c}\ar[d]& \mathbb{G}_m\ar[d]^-{z\mapsto z^d}\\
 \GL_{2d}\ar@{->>}[r]^-{\det}& \mathbb{G}_m
 }
 \]
 induces the commutative diagram
 \[
  \xymatrix{%
 X_*(T)\ar@{->>}[r]\ar[d]& X_*(\mathbb{G}_m)\ar[d]^-{\times d}\ar@{=}[r]& \pi_1(\GSp_{2d})\\
 X_*(T')\ar@{->>}[r]& X_*(\mathbb{G}_m)\ar@{=}[r]& \pi_1(\GL_{2d})\lefteqn{.}
 }
 \]
 In particular, the natural map $\pi_1(\GSp_{2d})\longrightarrow \pi_1(\GL_{2d})$ is injective.
\end{prf}

Therefore, there is a quasi-isogeny $\mathbb{X}\longrightarrow A[p^\infty]$ preserving polarizations.
If we replace $(\mathbb{X},\lambda_0)$ by the polarized $p$-divisible group 
$(A_0[p^\infty],\lambda_{A_0})$ associated with $(A_0,\lambda_{A_0})$,
the $G$-representation $H^i_{\mathrm{RZ}}$ remains unchanged. Thus, in order to prove Theorem \ref{thm:non-cusp-statement}, we may assume that $(\mathbb{X},\lambda_0)=(A_0[p^\infty],\lambda_{A_0})$.
In the remaining part of this article, we always assume it.
Moreover, we fix an isomorphism $H_1(A_0,\A^{\infty,p})\cong V'_{\A^{\infty,p}}$ preserving alternating pairings.

Denote the isogeny class of $(A_0,\lambda_{A_0})$ by $\phi$ and put $I^\phi=\Aut(\phi)$.
We have natural group homomorphisms $I^\phi\hooklongrightarrow J$ and
$I^\phi\hooklongrightarrow \Aut(H_1(A_0,\A^{\infty,p}))=\GSp(V'_{\A^{\infty,p}})$.
These are injective.

Let $Y_{K^p}$ be the reduced closed subscheme of $\overline{\Sh}_{K^p}$ such that 
$Y_{K^p}(\overline{\F}_p)$ consists of triples $(A,\lambda,\eta^p)$ where the $p$-divisible group
associated with $(A,\lambda)$ is isogenous to $(\X,\lambda_0)$.
It is the basic (or supersingular) stratum in the Newton stratification of $\overline{\Sh}_{K^p}$.
Note that $(A,\lambda,\eta^p)\in \overline{\Sh}_{K^p}(\overline{\F}_p)$ belongs to $Y_{K^p}(\overline{\F}_p)$
if and only if $(A,\lambda)\in \phi$ (\cite[Proposition 3.1.8]{MR2074714}, \cite[\S 7]{MR1124982}). 
We denote the formal completion of $\Sh_{K^p}$ along $Y_{K^p}$ by $(\Sh_{K^p})_{/Y_{K^p}}^\wedge$.

Now we can state the $p$-adic uniformization theorem:

\begin{thm}[{{\cite[Theorem 6.30]{MR1393439}}}]\label{thm:p-adic-unif}
 There exists a natural isomorphism of formal schemes:
 \[
  \theta_{K^p}\colon I^\phi\backslash (\M\times \GSp(V'_{\A^{\infty,p}})/K^p)\yrightarrow{\cong} (\Sh_{K^p})_{/Y_{K^p}}^\wedge.
 \]
 In the left hand side, $I^\phi$ acts on $\M$ through $I^\phi\hooklongrightarrow J$ and
 acts on $\GSp(V'_{\A^{\infty,p}})/K^p$ through $I^\phi\hooklongrightarrow \GSp(V'_{\A^{\infty,p}})$.
 
 The isomorphisms $\{\theta_{K^p}\}_{K^p}$ are compatible with change of $K^p$.
 (It is also compatible with the Hecke action of $\GSp_4(V'_{\A^{\infty,p}})$, but we do not use it.)
\end{thm}

Let us briefly recall the construction of the isomorphism $\theta_{K^p}$. 
Take a lift $(\widetilde{\mathbb{X}},\widetilde{\lambda}_0)$ of $(\mathbb{X},\lambda_0)$ over $\Z_{p^\infty}$
(such a lift is unique up to isomorphism).
Then, by the Serre-Tate theorem, the lift $(\widetilde{A}_0,\widetilde{\lambda}_{A_0})$ of $(A_0,\lambda_{A_0})$ is
canonically determined.
Let $S$ be an object of $\Nilp$, $(X,\rho)\in \M(S)$ and $[g]\in \GSp(V'_{\A^{\infty,p}})/K^p$.
Then $\rho$ extends uniquely to the quasi-isogeny
$\widetilde{\rho}\colon \widetilde{\mathbb{X}}\times_{\Z_{p^\infty}}S\longrightarrow X$.
We can see that there exist a polarized abelian variety $(A,\lambda)$
and a $p$-quasi-isogeny $\widetilde{A}_0\times_{\Z_{p^\infty}}S\longrightarrow A$ preserving polarizations,
such that the associated quasi-isogeny
$\widetilde{A}_0[p^\infty]\times_{\Z_{p^\infty}}S\longrightarrow A[p^\infty]$
coincides with $\widetilde{\rho}$. The fixed isomorphism 
$H_1(A_0,\A^{\infty,p})\cong V'_{\A^{\infty,p}}$ naturally induces a $K^p$-level structure $\eta$ of $A$.
The morphism $\theta_{K^p}$ is given by $\theta_{K^p}((X,\rho),[g])=(A,\lambda,\eta\circ g)$.

By composing the morphism $\M\longrightarrow \M\times \GSp(V'_{\A^{\infty,p}})/K^p$;
$x\longmapsto (x,[\id])$, we get a morphism
$\M\longrightarrow (\Sh_{K^p})_{/Y_{K^p}}^\wedge$, which is also denoted by $\theta_{K^p}$.
For $U\in \mathcal{Q}_0$, we denote the image of $U$ under $\theta_{K^p}$ by $Y_{K^p}(U)$.
It is an open subset of $Y_{K^p}$.

\begin{prop}\label{prop:p-adic-unif-M}
 Let $U$ be an element of $\mathcal{Q}_0$. Then for a sufficiently small compact open subgroup
 $K^p$ of $\GSp(V'_{\A^{\infty,p}})$, $\theta_{K^p}$ induces an isomorphism
 $U\yrightarrow{\cong} Y_{K^p}(U)$. Moreover, if we denote the open formal subscheme of
 $\M$ (resp.\ $(\Sh_{K^p})_{/Y_{K^p}}^\wedge$) whose underlying topological space is
 $U$ (resp.\ $Y_{K^p}(U)$) by $\M_{/U}$ (resp.\ $(\Sh_{K^p})_{/Y_{K^p}(U)}^\wedge)$,
 then $\theta_{K^p}$ induces an isomorphism 
 $\theta_{K^p}\colon \M_{/U}\yrightarrow{\cong} (\Sh_{K^p})_{/Y_{K^p}(U)}^\wedge$.
\end{prop}

\begin{prf}
 The proof is similar to \cite[Corollaire 3.1.4]{MR2074714}.
 Put $\Gamma_{K^p}=I^\phi\cap K^p$, where the intersection is taken in $\GSp(V'_{\A^{\infty,p}})$.
 It is known that $\Gamma_{K^p}$ is discrete and torsion-free \cite{MR1393439}.
 By Theorem \ref{thm:p-adic-unif}, $\theta_{K^p}$ gives an isomorphism from $\Gamma_{K^p}\backslash\M$
 to an open and closed formal subscheme of $(\Sh_{K^p})_{/Y_{K^p}}^\wedge$.
 By the same method as in \cite[Lemme 3.1.2, Proposition 3.1.3]{MR2074714}, we can see that 
 every element $\gamma\in \Gamma_{K^p}$ other than $1$ satisfies $\gamma\cdot U\cap U=\varnothing$
 if $K^p$ is sufficiently small. For such $K^p$, the natural morphism 
 $\M_{/U}\longrightarrow \Gamma_{K^p}\backslash\M$ is an open immersion.
 Thus we have an open immersion 
 $\M_{/U}\hooklongrightarrow \Gamma_{K^p}\backslash\M\yrightarrow[\cong]{\theta_{K^p}}(\Sh_{K^p})_{/Y_{K^p}}^\wedge$, whose image is $(\Sh_{K^p})_{/Y_{K^p}(U)}^\wedge$.
\end{prf}

Next we consider the case with Drinfeld level structures at $p$. 
Let $Y_{m,K^p}$ be the closed subscheme of $\overline{\Sh}_{m,K^p}$ obtained as the inverse image of $Y_{K^p}$
under $\overline{\Sh}_{m,K^p}\longrightarrow \overline{\Sh}_{K^p}$.
By the construction of $\theta_{K^p}$ described above, we have the following result:

\begin{cor}\label{cor:p-adic-unif-with-level}
 Let $m\ge 0$ be an integer.
 We can construct naturally a morphism 
 $\theta_{m,K^p}\colon \M_m\longrightarrow (\Sh_{m,K^p})^\wedge_{/Y_{m,K^p}}$ which makes the following diagram
 cartesian:
 \[
  \xymatrix{%
 \M_m\ar[r]^-{\theta_{m,K^p}}\ar[d]^-{p_{m0}}&(\Sh_{m,K^p})^\wedge_{/Y_{m,K^p}}\ar[d]\\
 \M\ar[r]^-{\theta_{K^p}}&(\Sh_{K^p})^\wedge_{/Y_{K^p}}\lefteqn{.}
 }
 \]
 In particular, the similar result as Proposition \ref{prop:p-adic-unif-M} holds for $\theta_{m,K^p}$;
 that is, for $U\in \mathcal{Q}_m$, $\theta_{m,K^p}$ induces $(\M_m)_{/U}\yrightarrow{\cong}(\Sh_{m,K^p})^\wedge_{/Y_{m,K^p}(U)}$
 if $K^p$ is sufficiently small.
\end{cor}

\section{Proof of the non-cuspidality result}\label{section:proof}
\subsection{The system of coefficients $\mathcal{F}^{[h]}$, $\mathcal{F}^{(h)}$}
\begin{defn}
 Let $m\ge 1$ and $0\le h\le 2$ be integers. We denote by $\mathcal{S}_{m,h}$ the set of
 direct summands of $L/p^mL$ of rank $4-h$, and  
 by $\mathcal{S}_{m,h}^{\mathrm{coi}}$ the subset of $\mathcal{S}_{m,h}$ consisting of
 coisotropic direct summands (recall that $I\in \mathcal{S}_{m,h}$ is said to be coisotropic if
 $I^\perp\subset I$).
 Put $\mathcal{S}_m=\bigcup_{h=0}^2\mathcal{S}_{m,h}$ and 
 $\mathcal{S}_m^{\mathrm{coi}}=\bigcup_{h=0}^2\mathcal{S}_{m,h}^{\mathrm{coi}}$.

 For $I\in \mathcal{S}_{m,h}$, let $\overline{\Sh}_{m,K^p,[I]}$ 
 be the $\overline{\F}_p$-scheme defined by 
 \begin{align*}
  \overline{\Sh}_{m,K^p,[I]}(S)&=\bigl\{(A,\lambda,\eta^p,\eta_p)\in \overline{\Sh}_{m,K^p,[I]}(S)
 \bigm| I\subset \Ker \eta_p\bigr\}.
 \end{align*}
 Clearly it is a closed subscheme of $\overline{\Sh}_{m,K^p}$.
 Similarly, we can define the closed formal subscheme $\M_{m,[I]}$ of $\M_m\otimes_{\Z_{p^\infty}}\overline{\F}_p$.
 Obviously, $\M_{m,[I]}$ is stable under the action of $J$ on $\M_m$.
\end{defn}

We denote by $Y_{m,K^p,[I]}$ the closed subscheme of $\overline{\Sh}_{m,K^p,[I]}$ obtained
as the inverse image of $Y_{m,K^p}$.
As Corollary \ref{cor:p-adic-unif-with-level}, we have the following cartesian diagram of formal schemes:
 \[
  \xymatrix{%
 \M_{m,[I]}\ar[r]\ar[d]&(\overline{\Sh}_{m,K^p,[I]})^\wedge_{/Y_{m,K^p,[I]}}\ar[d]\\
 \M_m\ar[r]^-{\theta_{m,K^p}}&(\Sh_{m,K^p})^\wedge_{/Y_{m,K^p}}\lefteqn{.}
 }
 \]

\begin{defn}
 For $I\in \mathcal{S}_m$,
 we put 
 \[
 \overline{\Sh}_{m,K^p,(I)}=\overline{\Sh}_{m,K^p,[I]}\setminus \bigcup_{I'\in \mathcal{S}_m,I\subsetneq I'} \overline{\Sh}_{m,K^p,[I']},
 \]
 which is an open subscheme of $\overline{\Sh}_{m,K^p,[I]}$, and thus is a subscheme of $\overline{\Sh}_{m,K^p}$.
 Moreover, for an integer $h$ with $0\le h\le 2$, we put
 $\overline{\Sh}^{[h]}_{m,K^p}=\bigcup_{I\in\mathcal{S}_{m,h}}\overline{\Sh}_{m,K^p,[I]}$
 and $\overline{\Sh}^{(h)}_{m,K^p}=\bigcup_{I\in\mathcal{S}_{m,h}}\overline{\Sh}_{m,K^p,(I)}$.
 The former is a closed subscheme of $\overline{\Sh}_{m,K^p}$, which is the scheme theoretic image of
 $\coprod_{I\in\mathcal{S}_{m,h}}\overline{\Sh}_{m,K^p,[I]}\longrightarrow \overline{\Sh}_{m,K^p}$.
 The latter is an open subscheme
 of $\overline{\Sh}^{[h]}_{m,K^p}$,
 since $\overline{\Sh}^{(h)}_{m,K^p}=\overline{\Sh}^{[h]}_{m,K^p}\setminus \overline{\Sh}^{[h-1]}_{m,K^p}$
 (if $h=0$, we put $\overline{\Sh}^{[-1]}_{m,K^p}=\varnothing$).
\end{defn}

\begin{lem}\label{lem:stratification-property}
 \begin{enumerate}
  \item Let $x=(A,\lambda,\eta^p,\eta_p)$ be an element of $\overline{\Sh}_{m,K^p}(\overline{\F}_p)$.
	Then, for $I\in \mathcal{S}_m$, $x\in \overline{\Sh}_{m,K^p,(I)}(\overline{\F}_p)$ if and
	only if $I=\Ker \eta_p$. For an integer $h$ with $0\le h\le 2$,
	$x\in \overline{\Sh}^{[h]}_{m,K^p}(\overline{\F}_p)$
	(resp.\ $x\in \overline{\Sh}^{(h)}_{m,K^p}(\overline{\F}_p)$) if and only if 
	$\rank_{\F_p} A[p]\le h$ (resp.\ $\rank_{\F_p} A[p]=h$).
  \item For every integer $h$ with $0\le h\le 2$,
	we have 
	$\overline{\Sh}^{(h)}_{m,K^p}=\coprod_{I\in \mathcal{S}_{m,h}}\overline{\Sh}_{m,K^p,(I)}$ as
	schemes.
  \item We have $(\overline{\Sh}^{[2]}_{m,K^p})_{\mathrm{red}}=(\overline{\Sh}_{m,K^p})_{\mathrm{red}}$ and
	$(\overline{\Sh}^{[0]}_{m,K^p})_{\mathrm{red}}=(Y_{m,K^p})_{\mathrm{red}}$.
 \end{enumerate}
\end{lem}

\begin{prf}
 Let us prove i).
 Put $X=A[p^\infty]$. Then there is an exact sequence $0\longrightarrow X_0\longrightarrow X\longrightarrow X_{\et}\longrightarrow 0$, where $X_0$ is a connected $p$-divisible group and $X_{\et}$ is an \'etale $p$-divisible group.
 By \cite[Lemma II.2.1]{MR1876802}, $\Ker \eta_p$ is a direct summand of $L/p^mL$ and
$(L/p^mL)/\Ker \eta_p\longrightarrow X_{\et}[p^m]$ is an isomorphism. 
 Thus $\Ker \eta_p\in \mathcal{S}_{m,r}$, where $r=\rank_{\Z/p^m\Z}X_{\et}[p^m]=\rank_{\F_p}A[p]\le 2$.
 By this, all the claims in i) are immediate.

 By i), $\overline{\Sh}^{(h)}_{m,K^p}$ coincides with
 $\coprod_{I\in \mathcal{S}_{m,h}}\overline{\Sh}_{m,K^p,(I)}$ as a set; thus to prove ii) it suffices to show that
 $\overline{\Sh}_{m,K^p,(I)}$ is closed (hence open) in $\overline{\Sh}^{(h)}_{m,K^p}$
 for every $I\in\mathcal{S}_{m,h}$. It is clear from 
 $\overline{\Sh}_{m,K^p,(I)}=\overline{\Sh}_{m,K^p,[I]}\cap \overline{\Sh}^{(h)}_{m,K^p}$.

 The former equality in iii) follows immediately from i). We will prove the latter.
 For $x=(A,\lambda,\eta^p,\eta_p)\in \overline{\Sh}^{[0]}_{m,K^p}(\overline{\F}_p)$,
 $X=A[p^\infty]$ has no \'etale part by i). Since $X^\vee\cong X$, $X$ has no multiplicative part. 
 Therefore $X$ is isoclinic of slope $1/2$; indeed, 
 if a Newton polygon with the terminal point $(4,2)$ has neither slope $0$ part nor slope $1$ part,
 then it is a line of slope $1/2$. Thus, by Lemma \ref{lem:RR},
 there is a quasi-isogeny $\mathbb{X}\longrightarrow X$ preserving polarizations; namely, 
 $x\in Y_{m,K^p}(\overline{\F}_p)$. The opposite inclusion is clear.
\end{prf}

\begin{rem}
 The latter part of iii) in Lemma \ref{lem:stratification-property}
 is the only place where the same argument does not work in the case $\GSp(2d)$ with $d\ge 3$.
\end{rem}

\begin{defn}
 Let $m\ge 1$ be an integer. Fix a compact open subgroup $K^p$ of $\GSp(V'_{\A^{\infty,p}})$.
 For $I\in \mathcal{S}_m$, denote the natural immersion
 $\overline{\Sh}_{m,K^p,(I)}\hooklongrightarrow \overline{\Sh}_{m,K^p}$ by $j_{m,I}$.
 For an integer $h$ with $0\le h\le 2$, denote the natural immersions
 $\overline{\Sh}_{m,K^p}^{[h]}\hooklongrightarrow \overline{\Sh}_{m,K^p}$ and
 $\overline{\Sh}_{m,K^p}^{(h)}\hooklongrightarrow \overline{\Sh}_{m,K^p}$
 by $j_m^{[h]}$ and $j_m^{(h)}$, respectively.
 
 We define $\mathcal{F}^\circ_{m,I}$, $\mathcal{F}_{m,I}$, $\mathcal{F}^{\circ [h]}_m$, $\mathcal{F}^{[h]}_m$,
 $\mathcal{F}^{\circ (h)}_m$ and $\mathcal{F}^{(h)}_m$ as follows:
 \begin{align*}
  \mathcal{F}^\circ_{m,I}&=\theta_{m}^*(Rj_{m,I*}Rj_{m,I}^!R\psi\Z_\ell)\vert_{Y_{m,K^p}},&
  \mathcal{F}_{m,I}&=\theta_{m}^*(Rj_{m,I*}Rj_{m,I}^!R\psi\Q_\ell)\vert_{Y_{m,K^p}},\\
  \mathcal{F}^{\circ [h]}_m&=\theta_{m}^*(Rj_{m*}^{[h]}Rj_m^{[h] !}R\psi\Z_\ell)\vert_{Y_{m,K^p}},&
  \mathcal{F}^{[h]}_m&=\theta_{m}^*(Rj_{m*}^{[h]}Rj_m^{[h] !}R\psi\Q_\ell)\vert_{Y_{m,K^p}},\\
  \mathcal{F}^{\circ (h)}_m&=\theta_{m}^*(Rj_{m*}^{(h)}Rj_m^{(h) !}R\psi\Z_\ell)\vert_{Y_{m,K^p}},&
  \mathcal{F}^{(h)}_m&=\theta_{m}^*(Rj_{m*}^{(h)}Rj_m^{(h) !}R\psi\Q_\ell)\vert_{Y_{m,K^p}}.
 \end{align*}
 Here $\theta_m\colon \Mbar_m\longrightarrow Y_{m,K^p}$ is the morphism
 induced from $\theta_{m,K^p}$ in Corollary \ref{cor:p-adic-unif-with-level}.

 These are independent of the choice of $K^p$; indeed, for another compact open subgroup $K'^p$
 contained in $K^p$, the natural map $\Sh_{m,K'^p}\longrightarrow \Sh_{m,K^p}$ is \'etale.
\end{defn}

\begin{prop}\label{prop:F-exact-seq}
 Let $h$ be an integer with $1\le h\le 2$.
 \begin{enumerate}
  \item We have the following distinguished triangle:
	\[
	 \mathcal{F}_{m}^{[h-1]}\longrightarrow \mathcal{F}_{m}^{[h]}\longrightarrow \mathcal{F}_{m}^{(h)}
	\longrightarrow \mathcal{F}_{m}^{[h-1]}[1].
	\]
  \item We have $\mathcal{F}^{(h)}_m=\bigoplus_{I\in\mathcal{S}_{m,h}}\mathcal{F}_{m,I}$.
 \end{enumerate}
 \end{prop}

\begin{prf}
 By the definition, i) is clear. ii) is also clear from Lemma \ref{lem:stratification-property} ii).
\end{prf}

\begin{prop}\label{prop:coisotropic}
 For $I\in \mathcal{S}_{m,h}\setminus\mathcal{S}_{m,h}^{\mathrm{coi}}$, we have
 $\mathcal{F}^\circ_{m,I}=\mathcal{F}_{m,I}=0$.
\end{prop}

\begin{prf}
 We will prove that $Rj_{m,I}^!R\psi\Z_\ell=0$. Since the dual of $Rj_{m,I}^!R\psi\Z_\ell$ is
 isomorphic to $j_{m,I}^*R\psi\Z_\ell(3)[6]$, it suffices to show that, for every
 $x\in \overline{\Sh}_{m,K^p,(I)}(\overline{\F}_p)$, no point on the generic fiber of $\Sh_{m,K^p}$
 specializes to $x$. In other words, for every complete discrete valuation ring $R$ with residue field 
 $\overline{\F}_p$ which is a flat $\Z_{p^\infty}$-algebra,
 and every $\Z_{p^\infty}$-morphism $\widetilde{x}\colon \Spec R\longrightarrow \Sh_{m,K^p}$, 
 the image of the closed point of $\Spec R$ under $\widetilde{x}$ does not lie in $\overline{\Sh}_{m,K^p,(I)}$.
 This is a consequence of the following lemma.
\end{prf}

\begin{lem}
 Let $R$ be a complete discrete valuation ring with perfect residue field $k$
 and with mixed characteristic $(0,p)$,
 and $(X,\lambda)$ a polarized $p$-divisible group over $R$. We denote the generic (resp.\ special) fiber of $X$
 by $X_\eta$ (resp.\ $X_s$). Then, for every $m\ge 1$, the kernel of the specialization map
 $X_\eta[p^m]\longrightarrow X_s[p^m]$ is a coisotropic direct summand of $X_\eta[p^m]$.
\end{lem}

\begin{prf}
 We shall prove that the kernel of the specialization map
 $T_pX_\eta\longrightarrow T_pX_s$ is a coisotropic direct summand of $T_pX_\eta$.
 Consider the exact sequence $0\longrightarrow X_{s,0}\longrightarrow X_s\longrightarrow X_{s,\et}\longrightarrow 0$ over $k$. It is canonically lifted to the exact sequence 
 $0\longrightarrow X_0\longrightarrow X\longrightarrow X_\et\longrightarrow 0$ over $R$,
 where $X_\et$ is an \'etale $p$-divisible group (\cf \cite[p.~76]{MR0347836}).
 Thus we have the following commutative diagram,
 whose rows are exact:
 \[
  \xymatrix{%
 0\ar[r]&T_pX_{0,\eta}\ar[r]\ar[d]&T_pX_{\eta}\ar[r]\ar[d]&T_pX_{\et,\eta}\ar[r]\ar[d]^-{\cong}&0\\
 0\ar[r]&0\ar[r]&T_pX_s\ar[r]&T_pX_{s,\et}\ar[r]&0\lefteqn{.}
 }
 \]
 Hence the kernel of $T_pX_\eta\longrightarrow T_pX_s$ coincides with $T_pX_{0,\eta}$.
 Therefore it suffices to show that the composite
 $(T_pX_{0,\eta})^\perp\longrightarrow T_pX_\eta\longrightarrow T_pX_{\et,\eta}$ is $0$.

 On the other hand, by the polarization $T_pX_\eta\yrightarrow{\cong} (T_pX_\eta)^\vee(1)$,
 $(T_pX_{0,\eta})^\perp$ corresponds to $(T_pX_{\et,\eta})^\vee(1)\cong T_pX_{\et,\eta}^\vee$.
 Thus it suffices to prove that every Galois-equivariant homomorphism
 $T_pX_{\et,\eta}^\vee\longrightarrow T_pX_{\et,\eta}$ is $0$.
 For this, we may replace the Tate modules $T_pX_{\et,\eta}^\vee$ and $T_pX_{\et,\eta}$ by 
 the rational Tate modules $V_pX_{\et,\eta}^\vee$ and $V_pX_{\et,\eta}$. These are crystalline representations
 and the corresponding filtered $\varphi$-modules are the rational Dieudonn\'e modules
 $D(X_{s,\et}^\vee)_{\Q}$ and $D(X_{s,\et})_{\Q}$, respectively. Since the slope of the former is $1$ and
 that of the latter is $0$, there is no $\varphi$-homomorphism other than $0$
 from $D(X_{s,\et}^\vee)_{\Q}$ to $D(X_{s,\et})_{\Q}$. This completes the proof.
\end{prf}

The following corollary is immediate from Proposition \ref{prop:F-exact-seq} ii) and Proposition \ref{prop:coisotropic}.

\begin{cor}\label{cor:coi-direct-sum}
 For $h$ with $1\le h\le 2$, we have $\mathcal{F}^{(h)}_m=\bigoplus_{I\in\mathcal{S}^{\mathrm{coi}}_{m,h}}\mathcal{F}_{m,I}$.
\end{cor}

Let us consider the action of $K_0$. Since $K_0/K_m$ naturally acts on $\Sh_{m,K^p}$ and
the action of $g\in K_0/K_m$ maps $\overline{\Sh}_{m,K^p,[I]}$ onto $\overline{\Sh}_{m,K^p,[g^{-1}I]}$,
the subschemes $\overline{\Sh}_{m,K^p}^{[h]}$ and $\overline{\Sh}_{m,K^p}^{(h)}$ are preserved
by the action of $K_0/K_m$. Therefore $\mathcal{F}_m^{\circ [h]}$, $\mathcal{F}_m^{[h]}$,
$\mathcal{F}_m^{\circ (h)}$ and $\mathcal{F}_m^{(h)}$ have natural $K_0/K_m$-equivariant structures.
Moreover, in the same way as in \cite[Proposition 2.5]{non-cusp},
we can observe that $\mathcal{F}^{[h]}=\{\mathcal{F}_m^{[h]}\}_{m\ge 1}$ and 
$\mathcal{F}^{(h)}=\{\mathcal{F}_m^{(h)}\}_{m\ge 1}$
(resp.\ $\mathcal{F}^{\circ [h]}=\{\mathcal{F}_m^{\circ [h]}\}_{m\ge 1}$ and 
$\mathcal{F}^{\circ (h)}=\{\mathcal{F}_m^{\circ (h)}\}_{m\ge 1}$) form systems of coefficients
(resp.\ integral coefficients) over $\{\Mbar_m\}_{m\ge 1}$.

Thanks to \cite{formalnearby}, we can define $J$-equivariant structures
on the systems of coefficients introduced above. 

\begin{prop}\label{prop:J-action}
 The complexes
 $\mathcal{F}^\circ_{m,I}$, $\mathcal{F}^{\circ [h]}_m$, $\mathcal{F}^{\circ (h)}_m$,
 $\mathcal{F}_{m,I}$, $\mathcal{F}^{[h]}_m$ and $\mathcal{F}^{(h)}_m$
 have natural $J$-equivariant structures.
 These structures are compatible with the distinguished triangles and the direct sum decompositions in
 Proposition \ref{prop:F-exact-seq}. 
\end{prop}

\begin{prf}
 We will prove the proposition for $\mathcal{F}_m^{(h)}$; other cases are similar.
 Put 
 \begin{align*}
  \overline{\Sh}_{m,K^p}^{[h]\wedge}&=(\Sh_{m,K^p})^\wedge_{/Y_{m,K^p}}\times_{\Sh_{m,K^p}}\overline{\Sh}_{m,K^p}^{[h]},&
 \overline{\Sh}_{m,K^p}^{(h)\wedge}&=\bigl(\overline{\Sh}_{m,K^p}^{[h]\wedge},\overline{\Sh}_{m,K^p}^{[h-1]\wedge}\bigr),\\
  \M_m^{[h]}&=\M_m\times_{(\Sh_{m,K^p})^\wedge_{/Y_{m,K^p}}}\overline{\Sh}_{m,K^p}^{[h]\wedge},&
  \M_m^{(h)}&=(\M_m^{[h]},\M_m^{[h-1]}).
 \end{align*}
 Then, by \cite[Proposition 3.11]{formalnearby}, we have the canonical isomorphism
 \[
  (Rj_{m*}^{(h)}Rj_m^{(h) !}R\psi\Q_\ell)\vert_{Y_{m,K^p}}\cong R\Psi_{\!(\Sh_{m,K^p})^\wedge_{/Y_{m,K^p}},\overline{\Sh}_{m,K^p}^{(h)\wedge}}\Q_\ell.
 \]
 Moreover, since $\theta_{m,K^p}$ is \'etale (\cf Corollary \ref{cor:p-adic-unif-with-level}), 
 by \cite[Proposition 3.14]{formalnearby}, we have the canonical isomorphism
 \[
 \mathcal{F}_m^{(h)}\cong R\Psi_{\!\!\!\M_m,\M_m^{(h)}}\Q_\ell.
 \]
 Since the action of $J$ on $\M_m$ preserves the closed formal subscheme $\M_{m,[I]}$ for every
 $I\in\mathcal{S}_m$, it also preserves the closed formal subscheme $\M_m^{[h]}$ for every $h$.
 Thus, by the functoriality \cite[Proposition 3.7]{formalnearby}, $R\Psi_{\!\!\!\M_m,\M_m^{(h)}}\Q_\ell$
 has a natural $J$-equivariant structure. We may import the structure into $\mathcal{F}^{(h)}_m$
 by the isomorphism above.

 The compatibilities with the exact sequence and the direct sum decomposition are clear
 from the construction (\cf \cite[Remark 3.8]{formalnearby}).
\end{prf}

It is easy to see that the actions defined in the previous proposition give $J$-equivariant structures
on the systems of (integral) coefficients $\mathcal{F}^{\circ [h]}$, $\mathcal{F}^{\circ (h)}$, 
$\mathcal{F}^{\circ [h]}$ and $\mathcal{F}^{\circ (h)}$. Thus we get the smooth representations
$H^i_c(\Mbar_\infty,\mathcal{F}^{[h]})$ and $H^i_c(\Mbar_\infty,\mathcal{F}^{(h)})$ of $K_0\times J$
(\cf Corollary \ref{cor:M_infty-coh}).

\begin{prop}\label{prop:RZ-F}
 There exists an isomorphism $H^i_c(\Mbar_\infty,\mathcal{F}^{[0]})\cong H^i_{\mathrm{RZ}}$,
 which is compatible with the action of $K_0\times J$.
\end{prop}

\begin{prf}
 Let $m\ge 1$ be an integer and $U\in \mathcal{Q}_m$.
 Then, by \cite[Corollary 4.40]{formalnearby} and Proposition \ref{prop:p-adic-unif-M},
 we have the $J$-equivariant isomorphism
 \[
  H^i_c(U,\mathcal{F}^{[0]}_m\vert_U)\cong H^i_c\bigl((\M_m)_{/U}^\rig\otimes_{\Q_{p^\infty}}\overline{\Q}_{p^\infty},\Q_\ell\bigr).
 \]
 Since this isomorphism is functorial, we have $K_0\times J$-equivariant isomorphisms
 \begin{gather*}
  H^i_c(\Mbar_m,\mathcal{F}_m^{[0]})\cong \varinjlim_{U\in\mathcal{Q}_m}H^i_c\bigl((\M_m)_{/U}^\rig\otimes_{\Q_{p^\infty}}\overline{\Q}_{p^\infty},\Q_\ell\bigr)
  \stackrel{(*)}{\cong} H^i_c(\M_m^\rig\otimes_{\Q_{p^\infty}}\overline{\Q}_{p^\infty},\Q_\ell),\\
  H^i_c(\Mbar_\infty,\mathcal{F}^{[0]})\cong \varinjlim_m H^i_c(\M_m^\rig\otimes_{\Q_{p^\infty}}\overline{\Q}_{p^\infty},\Q_\ell)=H^i_{\mathrm{RZ}}.
 \end{gather*}
 For the isomorphy of $(*)$, we need \cite[Proposition 2.1 (iv)]{MR1626021} and \cite[Lemma 4.14]{formalnearby}.
\end{prf}

\begin{rem}\label{rem:Berkovich-smoothness}
 We can deduce from Proposition \ref{prop:RZ-F} and Corollary \ref{cor:M_infty-coh}
 that the action of $K_0\times J$ on $H^i_{\mathrm{RZ}}$ is smooth.
\end{rem}

\subsection{$G$-action on $H^i_c(\Mbar_\infty,\mathcal{F}^{[h]})$, $H^i_c(\Mbar_\infty,\mathcal{F}^{(h)})$}
In this subsection, we define actions of $G$ on $H^i_c(\Mbar_\infty,\mathcal{F}^{[h]})$ and
$H^i_c(\Mbar_\infty,\mathcal{F}^{(h)})$ by using the method in \cite[\S 6]{MR2169874}.
Put $G^+=\{g\in G\mid g^{-1}L\subset L\}$, which is a submonoid of $G$.
For $g\in G^+$, let $e(g)$ be the minimal non-negative integer such that
$\Ker (g^{-1}\colon V/L\longrightarrow V/L)$ is contained in $p^{-e(g)}L/L$.
Since $\Ker g^{-1}=(gL+L)/L$, we have $gL\subset p^{-e(g)}L$.

In the sequel, we fix a compact open subgroup $K^p$ of $\GSp(V'_{\A^{\infty,p}})$ and denote
$\Sh_{m,K^p}$, $\overline{\Sh}_{m,K^p}$, $\overline{\Sh}_{m,K^p,[I]}$, \ldots by 
$\Sh_m$, $\overline{\Sh}_m$, $\overline{\Sh}_{m,[I]}$, \ldots, respectively.
Moreover, we fix $g\in G^+$ and denote $e(g)$ by $e$ for simplicity.

Assume that $m\ge e$. Let us consider the $\Z_{p^\infty}$-scheme $\Sh_{m,g}$ such that
for a $\Z_{p^\infty}$-scheme $S$, the set $\Sh_{m,g}(S)$ consists of isomorphism classes of 
quintuples $(A,\lambda,\eta^p,\eta_p,\mathcal{E})$ satisfying the following.
\begin{itemize}
 \item The quadruple $(A,\lambda,\eta^p,\eta_p)$ gives an element of $\Sh_m(S)$.
 \item $\mathcal{E}\subset X[p^e]$ is a finite flat subgroup scheme of order $p^{v_p(\det g^{-1})}$,
       where we put $X=A[p^\infty]$. It is self-dual with respect to $\lambda$, and satisfies
       $\eta'_p(\Ker g^{-1})\subset \mathcal{E}(S)$, where $\eta'_p$ denotes the composite
       $p^{-m}L/L\yrightarrow[\cong]{\times p^m}L/p^mL\yrightarrow{\eta_p}X[p^m]$.
 \item For $\mathcal{E}$ as above, we have the following commutative diagram:
       \[
	\xymatrix{%
       &p^{-m}L/L\ar[r]^-{\eta'_p}\ar@{->>}[d]_-{g^{-1}}& X[p^m]\ar@{^(-)}[r]\ar[d]& X\ar[d]\\
       &p^{-m}g^{-1}L/L\ar[r]& X[p^m]/\mathcal{E}\ar@{^(-)}[r]& X/\mathcal{E}\\
       L/p^{m-e}L\ar[r]^-{\cong}&p^{-m+e}L/L\ar[rr]\ar@{^(-)}[u]&& (X/\mathcal{E})[p^{m-e}]\ar@{^(-)}[u]\lefteqn{.}
       }
       \]
       We denote the composite of the lowest row by $\eta_p\circ g$ and assume that it gives a Drinfeld
       $(m-e)$-level structure.
\end{itemize}
We have the two natural morphisms 
\begin{gather*}
 \pr\colon \Sh_{m,g}\longrightarrow \Sh_m; (A,\lambda,\eta^p,\eta_p,\mathcal{E})\longmapsto (A,\lambda,\eta^p,\eta_p),\\
 [g]\colon \Sh_{m,g}\longrightarrow \Sh_{m-e}; (A,\lambda,\eta^p,\eta_p,\mathcal{E})\longmapsto (A/\mathcal{E},\lambda,\eta^p,\eta_p\circ g).
\end{gather*}
It is known that these are proper morphisms,
$\pr$ induces an isomorphism on the generic fibers, and 
$[g]$ induces the action of $g$ on the generic fibers \cite[Proposition 16, Proposition 17]{MR2169874}.

We can easily see that $\{\Sh_{m,g}\}_{m\ge e}$ form a projective system
whose transition maps are finite. Obviously, $\pr$ and $[g]$ are compatible with change of $m$. 

Similarly we can define the formal scheme $\M_{m,g}$ and the morphisms
$\pr\colon \M_{m,g}\longrightarrow \M_m$ and $[g]\colon \M_{m,g}\longrightarrow \M_{m-e}$.
The former morphism induces an isomorphism on the Raynaud generic fibers and 
the composite $[g]^\rig \circ (\pr^\rig)^{-1}$ coincides with the action of $g$.
The group $J$ naturally acts on $\M_{m,g}$ and two morphisms $\pr$ and $[g]$ are compatible
with the action of $J$.
Moreover, if we denote by $Y_{m,g}$ the inverse image of $Y_m\subset \Sh_m$ under $\pr\colon \Sh_{m,g}\longrightarrow \Sh_m$, then we can construct a morphism $\theta_{m,g}\colon \M_{m,g}\longrightarrow (\Sh_{m,g})^\wedge_{Y_{m,g}}$ which makes the following diagrams cartesian:
\[
  \xymatrix{%
 \M_{m,g}\ar[r]^-{\theta_{m,g}}\ar[d]^-{\pr}&(\Sh_{m,g})^\wedge_{/Y_{m,g}}\ar[d]^-{\pr}\\
 \M_m\ar[r]^-{\theta_m}&(\Sh_m)^\wedge_{/Y_m}\lefteqn{,}
 }\qquad
   \xymatrix{%
 \M_{m,g}\ar[r]^-{\theta_{m,g}}\ar[d]^-{[g]}&(\Sh_{m,g})^\wedge_{/Y_{m,g}}\ar[d]^-{[g]}\\
 \M_{m-e}\ar[r]^-{\theta_{m-e}}&(\Sh_{m-e})^\wedge_{/Y_{m-e}}\lefteqn{.}
 }
\]
Now let $h$ be an integer with $1\le h\le 2$ and $I\in \mathcal{S}_{m,h}$. Then we can define
the subschemes $\overline{\Sh}_{m,g,[I]}$, $\overline{\Sh}_{m,g,(I)}$, $\overline{\Sh}^{[h]}_{m,g}$ and
$\overline{\Sh}^{(h)}_{m,g}$ of $\Sh_{m,g}$
in the same way as $\overline{\Sh}_{m,[I]}$, $\overline{\Sh}_{m,(I)}$,
$\overline{\Sh}^{[h]}_m$ and $\overline{\Sh}^{(h)}_m$.
The following proposition is obvious:

\begin{prop}\label{prop:pr-strat}
 We have the commutative diagrams below:
\[
 \xymatrix{%
 \overline{\Sh}_{m,g,(I)}\ar[r]\ar[d]&\overline{\Sh}_{m,g,[I]}\ar[r]\ar[d]& \overline{\Sh}_{m,g}\ar[d]^-{\pr}\\
 \overline{\Sh}_{m,(I)}\ar[r]&\overline{\Sh}_{m,[I]}\ar[r]& \overline{\Sh}_m\lefteqn{,}
 }\qquad 
 \xymatrix{%
 \overline{\Sh}^{(h)}_{m,g}\ar[r]\ar[d]&\overline{\Sh}^{[h]}_{m,g}\ar[r]\ar[d]& \overline{\Sh}_{m,g}\ar[d]^-{\pr}\\
 \overline{\Sh}^{(h)}_{m}\ar[r]&\overline{\Sh}^{[h]}_m\ar[r]& \overline{\Sh}_m\lefteqn{.}
 }
\]
 The rectangles in the left diagram is cartesian. The rectangles in the right diagram is cartesian
 up to nilpotent elements (namely, $\overline{\Sh}_{m,g}^{[h]}\longrightarrow \overline{\Sh}_m^{[h]}\times_{\overline{\Sh}_m}\overline{\Sh}_{m,g}$ induces a homeomorphism on the underlying topological spaces, and so on).
\end{prop}

Let us consider how $\overline{\Sh}_{m,g,[I]}$ are mapped by $[g]\colon \Sh_{m,g}\longrightarrow \Sh_{m-e}$.
For this purpose, let us introduce some notation. 

\begin{defn}
 We denote by $\mathcal{S}_{\infty,h}$ the set of direct summands of $L$ of rank $4-h$ and
 by $\mathcal{S}_{\infty,h}^{\mathrm{coi}}$ the subset of $\mathcal{S}_{\infty,h}$
 consisting of coisotropic direct summands.
 We can identify $\mathcal{S}_{\infty,h}$ with the set of direct summands of $V$ of rank $4-h$;
 thus $G$ naturally acts on $\mathcal{S}_{\infty,h}$ and $\mathcal{S}^{\mathrm{coi}}_{\infty,h}$.
 Let $g^{-1}\colon \mathcal{S}_{m,h}\longrightarrow \mathcal{S}_{m-e,h}$ be the unique map which makes
 the following diagram commutative:
 \[
  \xymatrix{%
 \mathcal{S}_{\infty,h}\ar@{->>}[r]\ar[d]^-{g^{-1}}& \mathcal{S}_{m,h}\ar[d]^-{g^{-1}}\\
 \mathcal{S}_{\infty,h}\ar@{->>}[r]& \mathcal{S}_{m-e,h}\lefteqn{.}
 }
 \]
\end{defn}
The existence of such $g^{-1}$ follows from $p^mL\subset p^eL\subset g^{-1}L\subset L$.
Indeed, for direct summands $I$, $I'$ of $V$, we have
\begin{align*}
 &I\cap L+p^mL=I'\cap L+p^mL\implies g^{-1}I\cap g^{-1}L+p^mL=g^{-1}I'\cap g^{-1}L+p^mL\\
 &\qquad \implies g^{-1}I\cap g^{-1}L\cap p^eL+p^mL=g^{-1}I'\cap g^{-1}L\cap p^eL+p^mL\\
 &\qquad \iff g^{-1}I\cap p^eL+p^mL=g^{-1}I'\cap p^eL+p^mL\\
 &\qquad \iff g^{-1}I\cap L+p^{m-e}L=g^{-1}I'\cap L+p^{m-e}L.
\end{align*}
Obviously $g^{-1}\colon \mathcal{S}_{m,h}\longrightarrow \mathcal{S}_{m-e,h}$
induces a map from $\mathcal{S}_{m,h}^{\mathrm{coi}}$ to $\mathcal{S}_{m-e,h}^{\mathrm{coi}}$.

\begin{prop}\label{prop:g-strat}
 \begin{enumerate}
  \item For $h\in \{1,2\}$ and $I\in \mathcal{S}_{m,h}$, $[g]$ induces morphisms
	\begin{align*}
	 \Sh_{m,g,[I]}&\longrightarrow \Sh_{m-e,[g^{-1}I]}, &\Sh_{m,g,(I)}&\longrightarrow \Sh_{m-e,(g^{-1}I)},\\
	 \Sh_{m,g}^{[h]}&\longrightarrow \Sh_{m-e}^{[h]}, &\Sh_{m,g}^{(h)}&\longrightarrow \Sh_{m-e}^{(h)}.
	\end{align*}
  \item The rectangles of the following commutative diagram is cartesian up to nilpotent elements:
	\[
	  \xymatrix{%
	\overline{\Sh}^{(h)}_{m,g}\ar[r]\ar[d]&\overline{\Sh}^{[h]}_{m,g}\ar[r]\ar[d]& \overline{\Sh}_{m,g}\ar[d]^-{[g]}\\
	\overline{\Sh}^{(h)}_{m-e}\ar[r]&\overline{\Sh}^{[h]}_{m-e}\ar[r]& \overline{\Sh}_{m-e}\lefteqn{.}
	}
	\]
 \end{enumerate}
\end{prop}

\begin{prf}
 By the definition of $[g]$, it is clear that $[g]$ induces a morphism
 $\Sh_{m,g,[I]}\longrightarrow \Sh_{m-e,[g^{-1}I]}$ for $I\in \mathcal{S}_{m,h}$, and thus induces
 a morphism $\Sh_{m,g}^{[h]}\longrightarrow \Sh_{m-e}^{[h]}$.
 On the other hand, note that, for every $(A,\lambda,\eta^p,\eta_p,\mathcal{E})\in \Sh_{m,g}(\overline{\F}_p)$,
 the $p$-divisible groups $A[p^\infty]$ and $A[p^\infty]/\mathcal{E}$ are isogenous,
 and thus have the same \'etale heights. Therefore, by Lemma \ref{lem:stratification-property} i),
 the inverse image of $\overline{\Sh}_{m-e}^{[h]}$ (resp. $\overline{\Sh}_{m-e}^{(h)}$) under $[g]$
 coincides with $\overline{\Sh}_{m,g}^{[h]}$ (resp. $\overline{\Sh}_{m,g}^{(h)}$) as sets.
 Therefore a morphism $\Sh_{m,g}^{(h)}\longrightarrow \Sh_{m-e}^{(h)}$ is naturally induced and
 the rectangles in the diagram above are cartesian up to nilpotent elements.
 Finally, since $\overline{\Sh}_{m,g,(I)}=\overline{\Sh}_{m,g,[I]}\cap \overline{\Sh}_{m,g}^{(h)}$ and
 $\overline{\Sh}_{m-e,(g^{-1}I)}=\overline{\Sh}_{m-e,[g^{-1}I]}\cap \overline{\Sh}_{m-e}^{(h)}$, $[g]$ induces
 a morphism $\overline{\Sh}_{m,g,(I)}\longrightarrow \overline{\Sh}_{m-e,(g^{-1}I)}$.
\end{prf}

By Proposition \ref{prop:pr-strat} and Proposition \ref{prop:g-strat}, we have the natural cohomological
correspondence $\gamma_g$ from $\mathcal{F}^{[h]}_{m-e}$ (resp.\ $\mathcal{F}^{(h)}_{m-e}$) to 
$\mathcal{F}^{[h]}_m$ (resp.\ $\mathcal{F}^{(h)}_m$); see \S \ref{sec:coh-corr}.
This cohomological correspondence induces a homomorphism $\gamma_g$ from
 $H^i_c(\Mbar_{m-e},\mathcal{F}_{m-e}^{[h]})$
(resp.\ $H^i_c(\Mbar_{m-e},\mathcal{F}_{m-e}^{(h)})$) to $H^i_c(\Mbar_m,\mathcal{F}_m^{[h]})$
(resp.\ $H^i_c(\Mbar_m,\mathcal{F}_m^{(h)})$). Indeed, for $U\in \mathcal{Q}_{m-e}$,
we can take $U'\in\mathcal{Q}_m$ which contains $\pr([g]^{-1}(U))$. Then $\gamma_g$ induces 
$H^i_c(U,\mathcal{F}_{m-e}^{[h]}\vert_U)\longrightarrow H^i_c(U',\mathcal{F}_m^{[h]}\vert_{U'})$,
and therefore induces 
$H^i_c(\Mbar_{m-e},\mathcal{F}_{m-e}^{[h]})\longrightarrow H^i_c(\Mbar_m,\mathcal{F}_m^{[h]})$.
It is easy to see that this homomorphism is compatible with change of $m$;
hence we get the endomorphism $\gamma_g$ on $H^i_c(\Mbar_\infty,\mathcal{F}^{[h]})$
and $H^i_c(\Mbar_\infty,\mathcal{F}^{(h)})$. 

\begin{lem}\label{lem:G-J-comm}
 The endomorphism $\gamma_g$ commutes with the action of $J$ on $H^i_c(\Mbar_\infty,\mathcal{F}^{[h]})$
 and $H^i_c(\Mbar_\infty,\mathcal{F}^{(h)})$. 
\end{lem}

\begin{prf}
 We will only consider $\gamma_g$ on $H^i_c(\Mbar_\infty,\mathcal{F}^{[h]})$, since the other case is similar.
 Let $U\in \mathcal{Q}_{m-e}$ and $U'\in \mathcal{Q}_m$ be as above and put $W=[g]^{-1}(U)$, $W'=\pr^{-1}(U')$.
 It suffices to show the commutativity of the following diagram for $j\in J$:
 \[
  \xymatrix@!C=88pt{
  H^i_c(jU,\mathcal{F}^{[h]}_{m-e}\vert_{jU})\ar[r]^-{[g]^*}\ar[d]^-{j}&
  H^i_c(jW,\mathcal{F}^{[h]}_{m-e}\vert_{jW})\ar[r]\ar[d]^-{j}&
 H^i_c(jW',\mathcal{F}^{[h]}_{m-e}\vert_{jW'})\ar[r]^-{\pr_*}\ar[d]^-{j}&
 H^i_c(jU',\mathcal{F}^{[h]}_{m-e}\vert_{jU'})\ar[d]^-{j}\\
  H^i_c(U,\mathcal{F}^{[h]}_{m-e}\vert_U)\ar[r]^-{[g]^*}&
  H^i_c(W,\mathcal{F}^{[h]}_{m-e}\vert_W)\ar[r]&
  H^i_c(W',\mathcal{F}^{[h]}_{m-e}\vert_{W'})\ar[r]^-{\pr_*}&
  H^i_c(U',\mathcal{F}^{[h]}_{m-e}\vert_{U'})\lefteqn{.}
 }
 \]
 By the construction of the $J$-actions, the left and the middle rectangles are commutative.
 On the other hand, since $\pr$ is proper and induces an isomorphism on the generic fiber,
 $\pr_*$ is an isomorphism
 and its inverse is $\pr^*$. As $\pr^*$ commutes with the $J$-action, the right rectangle above 
 is also commutative. This concludes the proof.
\end{prf}

\begin{lem}\label{lem:action-of-G^+}
 \begin{enumerate}
  \item For $g,g'\in G^+$, $\gamma_{gg'}=\gamma_g\circ \gamma_{g'}$.
  \item For $g\in K_0$, $\gamma_g$ coincides with the action of $K_0$ on $H^i_c(\Mbar_\infty,\mathcal{F}^{[h]})$
	or $H^i_c(\Mbar_\infty,\mathcal{F}^{(h)})$, which we already introduced.
  \item The endomorphism $\gamma_{p^{-1}\cdot \id}$ an isomorphism (in fact, it coincides with the action
	of $p^{-1}\cdot \id\in J$).
 \end{enumerate}
\end{lem}

\begin{prf}
 i) follows from Corollary \ref{cor:compos-corr}. 
 ii) and iii) are consequences of \cite[Proposition 16, Proposition 17]{MR2169874} and the analogous
 properties for the Rapoport-Zink spaces (\cf \cite[Proposition 7.4 (4), (5)]{MR2074715}).
\end{prf}

Note that $G$ is generated by $G^+$ and $p\cdot \id$ as a monoid. Therefore, 
by the lemma above, we can extend the actions of $K_0$ on $H^i_c(\Mbar_\infty,\mathcal{F}^{[h]})$ and
$H^i_c(\Mbar_\infty,\mathcal{F}^{(h)})$ to whole $G$. Together with Lemma \ref{lem:G-J-comm}, we have
a smooth $G\times J$-module structures on $H^i_c(\Mbar_\infty,\mathcal{F}^{[h]})$ and
$H^i_c(\Mbar_\infty,\mathcal{F}^{(h)})$.
We can observe without difficulty that the isomorphism in Proposition \ref{prop:RZ-F} is in fact compatible
with the action of $G$:

\begin{prop}\label{prop:RZ-F-G}
 The isomorphism $H^i_c(\Mbar_\infty,\mathcal{F}^{[0]})\cong H^i_{\mathrm{RZ}}$ in Proposition \ref{prop:RZ-F}
 is an isomorphism of $G\times J$-modules.
\end{prop}

Next we investigate the $G$-module structure of $H^i_c(\Mbar_\infty,\mathcal{F}^{(h)})$ for $h\in \{1,2\}$.
Let us fix an element $\widetilde{I}(h)$ of $\mathcal{S}^{\mathrm{coi}}_{\infty,h}$ and denote its image
under the natural map $\mathcal{S}^{\mathrm{coi}}_{\infty,h}\longrightarrow \mathcal{S}^{\mathrm{coi}}_{m,h}$
by $\widetilde{I}(h)_m$. Put $P_h=\Stab_G(\widetilde{I}(h))$, which is a maximal parabolic subgroup of $G$.
Then we can identify $\mathcal{S}_{\infty,h}$ with $G/P_h=K_0/(P_h\cap K_0)$ and
$\mathcal{S}_{m,h}$ with $K_m\backslash G/P_h=K_m\backslash K_0/(P_h\cap K_0)$.
For $g\in G^+$ and an integer $m$ with $m\ge e:=e(g)$, 
$g^{-1}\colon \mathcal{S}_{m,h}\longrightarrow \mathcal{S}_{m-e,h}$ is identified with 
the map $K_m\backslash G/P_h\longrightarrow K_{m-e}\backslash G/P_h$; $K_mxP_h\longmapsto K_{m-e}g^{-1}xP_h$.

\begin{defn}
 We put $H^i_c(\Mbar_\infty,\mathcal{F}_{\widetilde{I}(h)})=\varinjlim_{m}H^i_c(\Mbar_m,\mathcal{F}_{m,\widetilde{I}(h)_m})$. Here the transition maps are given as follows: for integers $1\le m\le m'$,
 \begin{align*}
 H^i_c(\Mbar_m,\mathcal{F}_{m,\widetilde{I}(h)_m})&\longrightarrow H^i_c(\Mbar_{m'},p_{mm'}^*\mathcal{F}_{m,\widetilde{I}(h)_m})\longrightarrow \bigoplus_{\substack{I'\in \mathcal{S}^{\mathrm{coi}}_{m',h}\\ I'/p^mI'=\widetilde{I}(h)_m}}H^i_c(\Mbar_{m'},\mathcal{F}_{m',I'})\\
  &\twoheadlongrightarrow H^i_c(\Mbar_{m'},\mathcal{F}_{m',\widetilde{I}(h)_{m'}}).
 \end{align*}
\end{defn}
 It is easy to see that $H^i_c(\Mbar_\infty,\mathcal{F}_{\widetilde{I}(h)})$ has a structure of
 a smooth $P_h\times J$-module (use Theorem \ref{thm:J-fg-smooth} and Proposition \ref{prop:g-strat} i)).
 For each $m\ge 1$ we have the homomorphism
 \[
 H^i_c(\Mbar_m,\mathcal{F}^{(h)}_m)=\bigoplus_{I\in \mathcal{S}_{m,h}^{\mathrm{coi}}}H^i_c(\Mbar_m,\mathcal{F}_{m,I})\twoheadlongrightarrow H^i_c(\Mbar_m,\mathcal{F}_{m,\widetilde{I}(h)_m}),
 \]
 which induces the homomorphism 
 $H^i_c(\Mbar_\infty,\mathcal{F}^{(h)})\longrightarrow H^i_c(\Mbar_\infty,\mathcal{F}_{\widetilde{I}(h)})$.
 By Proposition \ref{prop:g-strat} i), we can prove that this is a homomorphism of $P_h\times J$-modules.

\begin{prop}\label{prop:induction-structure}
 We have an isomorphism 
 $H^i_c(\Mbar_\infty,\mathcal{F}^{(h)})\cong \Ind_{P_h}^GH^i_c(\Mbar_\infty,\mathcal{F}_{\widetilde{I}(h)})$ of
 $G\times J$-modules.
\end{prop}

\begin{prf}
 By the Frobenius reciprocity, we have a $G$-homomorphism
 $H^i_c(\Mbar_\infty,\mathcal{F}^{(h)})\longrightarrow \Ind_{P_h}^GH^i_c(\Mbar_\infty,\mathcal{F}_{\widetilde{I}(h)})$. We shall observe that this is bijective. For an integer $m\ge 1$, 
 we have 
 \begin{align*}
  H^i_c(\Mbar_m,\mathcal{F}_m^{(h)})&=\bigoplus_{I\in\mathcal{S}_{m,h}^{\mathrm{coi}}}H^i_c(\Mbar_m,\mathcal{F}_{m,I})=\bigoplus_{g\in K_m\backslash K_0/(P_h\cap K_0)}H^i_c(\Mbar_m,\mathcal{F}_{m,g^{-1}\widetilde{I}(h)_m})\\
  &\cong \Ind^{K_0/K_m}_{(P_h\cap K_0)/(P_h\cap K_m)}H^i_c(\Mbar_m,\mathcal{F}_{m,\widetilde{I}(h)_m}),
 \end{align*}
 where the last isomorphism, due to \cite[Lemme 13.2]{MR1719811}, is an isomorphism as $K_0$-modules.
 By taking the inductive limit, we have isomorphisms
 \[
  H^i_c(\Mbar_\infty,\mathcal{F}^{(h)})\yrightarrow{\cong} \Ind^{K_0}_{P_h\cap K_0}H^i_c(\Mbar_\infty,\mathcal{F}_{\widetilde{I}(h)_m})\yleftarrow{\cong} \Ind^{G}_{P_h}H^i_c(\Mbar_\infty,\mathcal{F}_{\widetilde{I}(h)_m})
 \]
 (the second isomorphy follows from the Iwasawa decomposition $G=P_hK_0$).
 By the proof of \cite[Lemme 13.2]{MR1719811}, it is easy to see that the first isomorphism above is
 nothing but the $K_0$-homomorphism obtained by the Frobenius reciprocity for $P_h\cap K_0\subset K_0$.
 Therefore the composite of the two isomorphisms above coincides with the $G$-homomorphism
 introduced at the beginning of this proof. Thus we conclude the proof. 
\end{prf}

\subsection{Proof of the main theorem}
We begin with the following result on non-cuspidality:

\begin{thm}\label{thm:ind-no-quasi-cuspidal}
 For every $i\in \Z$ and $h\in \{1,2\}$, the $G$-module $H^i_c(\Mbar_\infty,\mathcal{F}^{(h)})_{\overline{\Q}_\ell}$
 has no quasi-cuspidal subquotient.
\end{thm}

By Proposition \ref{prop:induction-structure} and \cite[2.4]{MR771671},
it suffices to show the following proposition:

\begin{prop}\label{prop:unip-rad-triv}
 Let $h\in \{1,2\}$.
 The unipotent radical $U_h$ of $P_h$ acts trivially on $H^i_c(\Mbar_\infty,\mathcal{F}_{\widetilde{I}(h)})_{\overline{\Q}_\ell}$.
\end{prop}

To prove Proposition \ref{prop:unip-rad-triv}, we need some preparations. 
In the sequel, let $\mathbf{G}$ and $\mathbf{H}$ be connected reductive groups over $\Q_p$,
$\mathbf{P}$ a parabolic subgroup of $\mathbf{G}$ and $\mathbf{U}$ the unipotent radical of $\mathbf{P}$.
We put $P=\mathbf{P}(\Q_p)$, $H=\mathbf{H}(\Q_p)$ and $U=\mathbf{U}(\Q_p)$.

\begin{lem}\label{lem:noether-adm}
 Let $A$ be a noetherian $\Q$-algebra and $V$ an $A$-module with a smooth $P$-action.
 Assume that $V$ is $A$-admissible in the sense of \cite[1.16]{MR771671}. Then $U$ acts on $V$ trivially.
\end{lem}

\begin{prf}
 First assume that $A$ is Artinian. Then we can prove the lemma in the same way as
 \cite[Lemme 13.2.3]{MR1719811} (we use length in place of dimension).

 For the general case, we use noetherian induction. Assume that the lemma holds for every proper quotient
 of $A$. Take a minimal prime ideal $\mathfrak{p}$ of $A$. Then $A_{\mathfrak{p}}$ is Artinian and
 $V_{\mathfrak{p}}$ is an $A_{\mathfrak{p}}$-admissible representation of $P$
 (note that $(V_{\mathfrak{p}})^K=(V^K)_{\mathfrak{p}}$ for every compact open subgroup $K$ of $P$).
 Therefore $U$ acts on $V_{\mathfrak{p}}$ trivially. Let $V'$ (resp.\ $V''$) be the kernel (resp.\ image)
 of $V\longrightarrow V_{\mathfrak{p}}$. Note that $V'$ and $V''$ are $A$-admissible representations of $P$,
 for $A$ is noetherian.

 Consider the following commutative diagram:
 \[
  \xymatrix{%
 0\ar[r]& V'\ar[r]\ar[d]^-{(1)}& V\ar[r]\ar[d]^-{(2)}& V''\ar[r]\ar[d]^-{(3)}& 0\\
 0\ar[r]& V'_U\ar[r]& V_U\ar[r]& (V_{\mathfrak{p}})_U\lefteqn{.}&
 }
 \]
 It is well-known that the functor taking $U$-coinvariant $V\longmapsto V_U$ is an exact functor;
 thus the lower row in the diagram above is exact. On the other hand, the arrow labeled $(3)$ is injective,
 since it is the composite of $V''\hooklongrightarrow V_{\mathfrak{p}}\yrightarrow{\cong}(V_{\mathfrak{p}})_U$.
 Therefore, by the snake lemma, the injectivity of $(2)$ is equivalent to that of $(1)$.
 In other words, we have only to prove that the action of $U$ on $V'$ is trivial.

 On the other hand, by the definition, $V'$ is the union of $V_s:=\{x\in V\mid sx=0\}$
 for $s\in A\setminus \mathfrak{p}$. Since $V_s$ can be regarded as an admissible $A/(s)$-representation,
 $U$ acts on $V_s$ trivially by the induction hypothesis. Hence $U$ acts on $V'$ trivially.
\end{prf}

\begin{prop}\label{Prop:Bernstein-center}
 Let $V$ be a smooth representation of $P\times H$ over $\overline{\Q}_\ell$ and assume that
 for every compact open subgroup $K$ of $P$, $V^K$ is a finitely generated $H$-module.
 Then $U$ acts on $V$ trivially.
\end{prop}

\begin{prf}
 Since $\overline{\Q}_\ell$ and $\C$ are isomorphic as fields, we may replace $\overline{\Q}_\ell$ 
 in the statement by $\C$. Let $\mathfrak{Z}$ be the Bernstein center of $H$ \cite{MR771671}. 
 It is decomposed as $\mathfrak{Z}=\prod_{\theta\in\Theta}\mathfrak{Z}_\theta$,
 where $\Theta$ denotes the set of connected components of the Bernstein variety of $H$.
 For $\theta\in \Theta$, we denote the $\theta$-part of $V$ by $V_\theta$.
 Then we have the canonical decomposition $V=\bigoplus_{\theta\in\Theta}V_\theta$, which is
 compatible with the action of $P\times H$.
 Therefore, by replacing $V$ with $V_\theta$, 
 we may assume that the action of $\mathfrak{Z}$ on $V$ factors through
 $\mathfrak{Z}_\theta$ for some $\theta\in \Theta$.

 By the assumption and \cite[Proposition 3.3]{MR771671}, for every compact open subgroup $K$ of $P$,
 $V^K$ is a $\mathfrak{Z}_\theta$-admissible $H$-module.
 Namely, for every compact open subgroup $K$ (resp.\ $K'$) of $P$ (resp.\ $H$),
 $V^{K\times K'}$ is a finitely generated $\mathfrak{Z}_\theta$-module.
 In other words, for every compact open subgroup $K'$ of $H$, $V^{K'}$ is a $\mathfrak{Z}_\theta$-admissible
 $P$-module. Since $\mathfrak{Z}_\theta$ is a finitely generated $\C$-algebra, $U$ acts trivially 
 on $V^{K'}$ by Lemma \ref{lem:noether-adm}. Therefore $U$ acts trivially on $V$ also.
\end{prf}

\begin{prf}[of Proposition \ref{prop:unip-rad-triv}]
 By Proposition \ref{Prop:Bernstein-center}, we have only to prove that, for every $m\ge 1$,
 $H^i_c(\Mbar_\infty,\mathcal{F}_{\widetilde{I}(h)})^{P_h\cap K_m}$ is a finitely generated $J$-module
 (recall that a finitely generated $J$-module is noetherian \cite[Remarque 3.12]{MR771671}).
 As a $J$-module, it is a direct summand of
 $(\Ind_{P_h}^G H^i_c(\Mbar_\infty,\mathcal{F}_{\widetilde{I}(h)}))^{K_m}\cong H^i_c(\Mbar_\infty,\mathcal{F}^{(h)})^{K_m}$. On the other hand, by Corollary \ref{cor:M_infty-coh}, 
 $H^i_c(\Mbar_\infty,\mathcal{F}^{(h)})^{K_m}$ is a finitely generated $J$-module. Thus
 $H^i_c(\Mbar_\infty,\mathcal{F}_{\widetilde{I}(h)})^{P_h\cap K_m}$ is also finitely generated.
\end{prf}

\begin{prop}\label{prop:coh-vanishing}
 Let $i$ be an integer. If $i\ge 5$, then $H^i_c(\Mbar_{\infty},\mathcal{F}^{[2]})=0$.
 On the other hand, if $i\le 1$, then $H^i_c(\Mbar_{\infty},\mathcal{F}^{[0]})=0$.
\end{prop}

\begin{prf}
 By the definition, it suffices to show that for every $m\ge 1$ and every $U\in\mathcal{Q}_m$
 we have $H^i_c(U,\mathcal{F}^{[2]}_m\vert_U)=0$ for $i\ge 5$ and 
 $H^i_c(U,\mathcal{F}^{[0]}_m\vert_U)=0$ for $i\le 1$. Thus the claim is reduced to the following
 lemma.
\end{prf}

\begin{lem}
 Let $S$ be the spectrum of a strict henselian discrete valuation ring and $X$ a separated $S$-scheme
 of finite type. We denote its special (resp.\ generic) fiber by
 $X_s$ (resp.\ $X_\eta$). Let $Z$ be a closed subscheme of $X_s$ and denote the natural closed
 immersion $Z\hooklongrightarrow X$ by $i$. Assume that $X_\eta$ is smooth of pure dimension $d$
 and $Z$ is purely $d'$-dimensional.

 Then we have $H^n(Z,i^*R\psi_X\Q_\ell)=H^n_c(Z,i^*R\psi_X\Q_\ell)=0$ for $n>d+d'$ and 
 $H^n_c(Z,Ri^!R\psi_X\Q_\ell)=0$ for $n<d-d'$.
\end{lem}

\begin{prf}
 First note that $H^n(Z,i^*R^k\psi_X\Q_\ell)=H_c^n(Z,i^*R^k\psi_X\Q_\ell)=0$ if $n>2\dim Z$ 
 or $n>2\dim (\supp R^k\psi_X\Q_\ell)$.
 By \cite[Proposition 4.4.2]{MR751966}, for each $k\ge 0$ we have $\dim(\supp R^k\psi_X\Q_\ell)\le d-k$;
 therefore if $n+k>d+d'$ then we have
 \begin{align*}
  n&>d'+(d-k)\ge \dim Z+\dim(\supp R^k\psi_X\Q_\ell)\\
  &\ge 2\min \{\dim Z,\dim(\supp R^k\psi_X\Q_\ell)\}
 \end{align*}
 and thus $H^n(Z,i^*R^k\psi_X\Q_\ell)=H_c^n(Z,i^*R^k\psi_X\Q_\ell)=0$.
 By the spectral sequence,
 we have $H^n(Z,i^*R\psi_X\Q_\ell)=H^n_c(Z,i^*R\psi_X\Q_\ell)=0$ for $n>d+d'$.

 On the other hand, by the Poincar\'e duality, we have
 \begin{align*}
  H^n_c(Z,Ri^!R\psi_X\Q_\ell)&=H^{-n}(Z,D_Z(Ri^!R\psi_X\Q_\ell))^{\vee}
  =H^{-n}(Z,i^*R\psi_XD_{X_\eta}\Q_\ell)^{\vee}\\
  &=H^{-n}(Z,i^*R\psi_X\Q_\ell(d)[2d])^{\vee}
  =H^{2d-n}(Z,i^*R\psi_X\Q_\ell)^{\vee}(-d),
 \end{align*}
 where $D_Z$ (resp.\ $D_{X_\eta}$) denotes the dualizing functor with respect to $Z$ (resp.\ $X_\eta$).
 Therefore it vanishes if $2d-n>d+d'$, namely, $n<d-d'$.
\end{prf}

Now we can prove our main theorem.

\begin{prf}[of Theorem \ref{thm:non-cusp-statement}]
 By Proposition \ref{prop:RZ-F} and Proposition \ref{prop:coh-vanishing}, we have
 $H^i_{\mathrm{RZ}}=0$ for $i\le 1$. Therefore we may assume that $i\ge 5$.

 By Proposition \ref{prop:F-exact-seq} i), we have the exact sequence of smooth $G$-modules
 \[
 H_c^{i-1}(\Mbar_\infty,\mathcal{F}^{(h)})_{\overline{\Q}_\ell}\longrightarrow H_c^i(\Mbar_\infty,\mathcal{F}^{[h-1]})_{\overline{\Q}_\ell}\longrightarrow H^i_c(\Mbar_\infty,\mathcal{F}^{[h]})_{\overline{\Q}_\ell}
 \]
 for every $h$ with $1\le h\le 2$. Moreover, $H_c^i(\Mbar_\infty,\mathcal{F}^{(h)})_{\overline{\Q}_\ell}$ has
 no quasi-cuspidal subquotient by Theorem \ref{thm:ind-no-quasi-cuspidal}.
 Thus, starting from $H^i_c(\Mbar_{\infty},\mathcal{F}^{[2]})_{\overline{\Q}_\ell}=0$
 (Proposition \ref{prop:coh-vanishing}), we can inductively prove that 
 $H^i_c(\Mbar_{\infty},\mathcal{F}^{[h]})_{\overline{\Q}_\ell}$ has no quasi-cuspidal subquotient;
 indeed, the property that a representation has no quasi-cuspidal subquotient is
 stable under sub, quotient and extension (use the canonical decomposition in \cite[2.3.1]{MR771671}).
 In particular, $H^i_c(\Mbar_\infty,\mathcal{F}^{[0]})_{\overline{\Q}_\ell}\cong H^i_{\mathrm{RZ},\overline{\Q}_\ell}$ (\cf Proposition \ref{prop:RZ-F-G})
 has no quasi-cuspidal subquotient. This completes the proof.
\end{prf}

\section{Appendix: Complements on cohomological correspondences}\label{sec:coh-corr}
In this section, we recall the notion of cohomological correspondences (\cf \cite[Expos\'e III]{SGA5}, \cite{MR1431137})
and give some results on them. These are used to define the action of $G$
on $H^i_c(\Mbar_\infty,\mathcal{F}^{[h]})$ and $H^i_c(\Mbar_\infty,\mathcal{F}^{(h)})$.

In this section, we change our notation.
Let $k$ be a field and $\ell$ a prime number which is invertible in $k$.
We denote one of $\Z/\ell^n\Z$ or $\Q_\ell$ by $\Lambda$.
Let $X_1$ and $X_2$ be schemes which are separated of finite type over $k$,
and $L_i$ an object of $D^b_c(X_i,\Lambda)$ for $i=1,2$ respectively.
A \textit{cohomological correspondence} from $L_1$ to $L_2$ is a pair $(\gamma,c)$
consisting of a separated $k$-morphism of finite type $\gamma\colon \Gamma\longrightarrow X_1\times X_2$ and
a morphism $c\colon \gamma_1^*L_1\longrightarrow R\gamma_2^!L_2$ in the category $D_c^b(\Gamma,\Lambda)$,
where we denote $\pr_i\circ \gamma$ by $\gamma_i$. 
For simplicity, we also write $c$ for $(\gamma,c)$, if there is no risk of confusion.
If we are given a cohomological correspondence $(\gamma,c)$ where $\gamma_1$ is proper,
then we have the associated morphism $R\Gamma_c(c)\colon R\Gamma_c(X_1,L_1)\longrightarrow R\Gamma_c(X_2,L_2)$
by composing
\begin{align*}
 R\Gamma_c(X_1,L_1)&\yrightarrow{\gamma_1^*}R\Gamma_c(\Gamma,\gamma_1^*L_1)\yrightarrow{R\Gamma_c(c)}
 R\Gamma_c(\Gamma,R\gamma_2^!L_1)=R\Gamma_c(X_2,R{\gamma_2}_!R\gamma_2^!L_2)\\
 &\yrightarrow{\adj}R\Gamma_c(X_2,L_2).
\end{align*}

We can compose two cohomological correspondences. Let $X_3$ be another scheme
which is separated of finite type over $k$ and $L_3\in D^b_c(X_3,\Lambda)$.
Let $(\gamma',c')$ be a cohomological correspondence from $L_2$ to $L_3$.
Consider the following diagram
 \[
  \xymatrix{%
  \Gamma\times_{X_2}\Gamma'\ar[r]^-{\pr_2}\ar[d]^-{\pr_1}&\Gamma'\ar[r]^-{\gamma'_2}\ar[d]^-{\gamma'_1}&X_3\\
  \Gamma\ar[r]^-{\gamma_2}\ar[d]^-{\gamma_1}& X_2&\\
  X_1\lefteqn{.}
  }
 \]
Let $\gamma''$ be the natural morphism $\Gamma\times_{X_2}\Gamma'\longrightarrow X_1\times X_3$
and $c''\colon \gamma_1''^*L_1\longrightarrow R\gamma_2''^!L_3$ the composite of
\[
 \gamma_1''^*L_1=\pr_1^*\gamma_1^*L_1\yrightarrow{\pr_1^*(c)}\pr_1^*R\gamma_2^!L_2\yrightarrow{\mathrm{b.c.}}
 R\pr_2^!\gamma'^*_1L_2\yrightarrow{R\pr_2^!(c')}R\pr_2^!R\gamma'^!_2L_3=R\gamma_2''^!L_3,
\]
where $\mathrm{b.c.}$ denotes the base change morphism.
We call the cohomological correspondence $(\gamma'',c'')$ the composite of $(\gamma,c)$ and $(\gamma',c')$,
and denote it by $c'\circ c$. It is not difficult to see that if $\gamma_1$ and $\gamma'_1$ are proper, then 
$\gamma_1''$ is also proper and $R\Gamma_c(c'\circ c)=R\Gamma_c(c')\circ R\Gamma_c(c)$.

Let us recall some operations for cohomological correspondences.
Let $X_1$, $X_2$, $X_1'$ and $X_2'$ be schemes which are separated of finite type over $k$,
and $\gamma\colon \Gamma\longrightarrow X_1\times X_2$ and $\gamma'\colon \Gamma'\longrightarrow X'_1\times X'_2$
separated $k$-morphisms of finite type. Assume that the following commutative diagram is given:
\[
 \xymatrix{%
 X_1'\ar[d]^-{a_1}& \Gamma'\ar[d]^-{a}\ar[l]_-{\gamma_1'}\ar[r]^-{\gamma_2'}& X_2'\ar[d]^-{a_2}\\
 X_1& \Gamma\ar[l]_-{\gamma_1}\ar[r]^-{\gamma_2}& X_2\lefteqn{.}
 }
\]
First assume that every vertical morphism is proper.
Let $L'_i$ be an object of $D^b_c(X_i',\Lambda)$
for each $i=1,2$ and $(\gamma',c')$ a cohomological correspondence from $L_1'$ to $L_2'$.
Then we can define the cohomological correspondence $(\gamma,a_*c')$ from $R{a_1}_*L_1'$ to $R{a_2}_*L_2'$ by
\begin{align*}
 \gamma_1^*R{a_1}_*L_1'&\yrightarrow{\mathrm{b.c.}}Ra_*\gamma'^*_1L_1'\yrightarrow{Ra_*(c')}
 Ra_*R\gamma'^!_2L_2'=Ra_!R\gamma'^!_2L_2'\\
 &\yrightarrow{\mathrm{b.c.}} R\gamma^!_2R{a_2}_!L_2'=R\gamma^!_2R{a_2}_*L_2'.
\end{align*}
The cohomological correspondence $(\gamma,a_*c')$ is called the push-forward of $(\gamma',c')$ by $a$.
It is easy to see that push-forward is compatible with composition.
Moreover, we have the following lemma whose proof is also immediate:

\begin{lem}\label{lem:coh-corr-push}
 In the above diagram, assume that $X_1=X_1'$, $X_2=X_2'$, $a_1=a_2=\id$ and $a$ is proper.
 Let $L_i$ be an object of $D^b_c(X_i,\Lambda)$ for each $i=1,2$ and $(\gamma',c')$ a cohomological
 correspondence from $L_1$ to $L_2$. Then we have $R\Gamma_c(a_*c')=R\Gamma_c(c')$.
\end{lem}

Next we assume that the right rectangle in the diagram above is cartesian.
Let $L_i'$ and $(\gamma',c')$ be as above. Then we have the cohomological correspondence
$(\gamma,a_*c')$ from $R{a_1}_*L_1'$ to $R{a_2}_*L_2'$ by 
\[
 \gamma_1^*R{a_1}_*L_1'\yrightarrow{\mathrm{b.c.}}Ra_*\gamma'^*_1L_1'\yrightarrow{Ra_*(c')}
 Ra_*R\gamma'^!_2L_2'\yrightarrow{\mathrm{b.c.}} R\gamma^!_2R{a_2}_*L_2'.
\]
On the other hand, let $L_i$ be an object of $D^b_c(X_i,\Lambda)$ for each $i=1,2$ and $(\gamma,c)$
a cohomological correspondence from $L_1$ to $L_2$.
Then we have the cohomological correspondence
$(\gamma,a^*c)$ from $a_1^*L_1$ to $a_2^*L_2$ by 
\[
 \gamma_1'^*a_1^*L_1=a^*\gamma^*_1L_1\yrightarrow{a^*(c)}
 a^*R\gamma^!_2L_2\yrightarrow{\mathrm{b.c.}} R\gamma'^!_2a_2^*L_2.
\]

Finally assume that the left rectangle in the diagram above is cartesian.
Let $L_i$ and $(\gamma,c)$ be as above. 
Then we have the cohomological correspondence
$(\gamma,a^!c)$ from $Ra_1^!L_1$ to $Ra_2^!L_2$ by 
\[
 \gamma_1'^*Ra_1^!L_1\yrightarrow{\mathrm{b.c.}}Ra^!\gamma^*_1L_1\yrightarrow{Ra^!(c)}
 Ra^!R\gamma^!_2L_2=R\gamma'^!_2Ra_2^!L_2.
\]
These constructions are also compatible with composition.

Next we recall the specialization of cohomological correspondences.
Let $S$ be the spectrum of a strict henselian discrete valuation ring on which $\ell$ is invertible.
For an $S$-scheme $X$, we denote its special (resp.\ generic) fiber by
$X_s$ (resp.\ $X_{\eta}$).

Let $X_1$, $X_2$ be schemes which are separated of finite type over $S$ and
$\gamma\colon \Gamma\longrightarrow X_1\times_S X_2$ a separated $S$-morphism of finite type.
Let $L_i$ be an object of $D^b_c(X_{i,{\eta}},\Lambda)$ for each $i=1,2$
and $(\gamma_{\eta},c)$ a cohomological correspondence from $L_1$ to $L_2$.
Then we have the cohomological correspondence $(\gamma_s,R\psi(c))$ from
$R\psi L_1$ to $R\psi L_2$ by
\[
 \gamma_{1,s}^*R\psi L_1\longrightarrow R\psi \gamma_{1,\eta}^*L_1\yrightarrow{R\psi(c)}
 R\psi R\gamma_{2,\eta}^!L_2 \longrightarrow R\gamma_{2,s}^!R\psi L_2.
\]
It is easy to see that this construction is compatible with composition and proper push-forward
(\cf \cite[Proposition 1.6.1]{MR1431137}).

Now we will give the main result in this section. Let $X_i$, $\gamma$, $L_i$ be as above and 
$Y_i$ (resp.\ $Z_i$) a closed (resp.\ locally closed) subscheme of $X_{i,s}$.
Assume that $\gamma_{1,s}^{-1}(Y_1)=\gamma_{2,s}^{-1}(Y_2)$ and $\gamma_{1,s}^{-1}(Z_1)=\gamma_{2,s}^{-1}(Z_2)$
as subschemes of $\Gamma_s$, and denote the former by $\Gamma_Y$ and the latter by $\Gamma_Z$.
Then we have the following diagrams whose rectangles are cartesian:
\[
 \xymatrix{%
 Y_1\ar[d]^-{i_1}& \Gamma_Y\ar[d]^-{i}\ar[l]_-{\gamma_{Y,1}}\ar[r]^-{\gamma_{Y,2}}& Y_2\ar[d]^-{i_2}\\
 X_{1,s}& \Gamma_s\ar[l]_-{\gamma_{1,s}}\ar[r]^-{\gamma_{2,s}}& X_{2,s}\lefteqn{,}
 }\qquad\qquad
 \xymatrix{%
 Z_1\ar[d]^-{j_1}& \Gamma_Z\ar[d]^-{j}\ar[l]_-{\gamma_{Z,1}}\ar[r]^-{\gamma_{Z,2}}& Z_2\ar[d]^-{j_2}\\
 X_{1,s}& \Gamma_s\ar[l]_-{\gamma_{1,s}}\ar[r]^-{\gamma_{2,s}}& X_{2,s}\lefteqn{.}
 }
\]
Therefore, for a cohomological correspondence $(\gamma_{\eta},c)$ from $L_1$ to $L_2$,
the cohomological correspondence $i^*j_*j^!R\psi(c)$ from $i_1^*R{j_1}_*Rj_1^!R\psi L_1$ to
 $i_2^*R{j_2}_*Rj_2^!R\psi L_2$ is induced. If moreover we assume that $\gamma_1$ is proper, then we have
\[
 R\Gamma_c(i^*j_*j^!R\psi(c))\colon R\Gamma_c(X_{1,s},i_1^*R{j_1}_*Rj_1^!R\psi L_1)\longrightarrow R\Gamma_c(X_{2,s},i_2^*R{j_2}_*Rj_2^!R\psi L_2).
\]

\begin{prop}\label{prop:specialization-indep}
 The morphism $R\Gamma_c(i_*j_*j^!R\psi(c))$ depends only on the cohomological correspondence
 $(\gamma_{\eta},c)$. More precisely, if another $S$-morphism 
 $\gamma'\colon \Gamma'\longrightarrow X_1\times_SX_2$ has the same generic fiber as $\gamma$
 and satisfies the conditions that $\gamma'^{-1}_{1,s}(Y_1)=\gamma'^{-1}_{2,s}(Y_2)$,
 $\gamma'^{-1}_{1,s}(Z_1)=\gamma'^{-1}_{2,s}(Z_2)$ and $\gamma'_1$ is proper, then 
 the morphism $R\Gamma_c(i'^*j'_*j'^!R\psi(c))$ induced from $\gamma'$ is equal to
 $R\Gamma_c(i^*j_*j^!R\psi(c))$ (here $i'$ and $j'$ are defined in the same way as $i$ and $j$).
\end{prop}

\begin{prf}
 Since $\Gamma$ and $\Gamma'$ have the same generic fiber, there is the ``diagonal''
 in the generic fiber of $\Gamma\times_{X_1\times_SX_2}\Gamma'$. Let $\Gamma''$ be the closure
 of it in $\Gamma\times_{X_1\times_SX_2}\Gamma'$.
 Then $\Gamma''$ has the same generic fiber as $\Gamma$.
 We have natural morphisms $\Gamma''\longrightarrow \Gamma$ and $\Gamma''\longrightarrow \Gamma'$, which are
 proper since $\gamma$ and $\gamma'$ are proper. 
 Therefore $\gamma''\colon \Gamma''\longrightarrow X_1\times_SX_2$ also satisfies the same conditions as
 $\gamma$ and $\gamma'$. By replacing $\gamma'$ by $\gamma''$, we may assume that there exists
 a proper morphism $a\colon \Gamma'\longrightarrow \Gamma$ such that $\gamma\circ a=\gamma'$.

 Then, it is easy to see that the push-forward of the cohomological correspondence
 $(\gamma_s',i'^*j'_*j'^!R\psi(c))$ by $a_s$ coincides with $(\gamma_s,i^*j_*j^!R\psi(c))$.
 Therefore the proposition follows from Lemma \ref{lem:coh-corr-push}.
\end{prf}

\begin{cor}\label{cor:compos-corr}
 Let $X_1$, $X_2$ and $X_3$ be schemes which are separated of finite type over $S$,
 $Y_i$ (resp.\ $Z_i$) a closed (resp.\ locally closed) subscheme of $X_i$,
 and $L_i$ an object of $D^b_c(X_{i,\eta},\Lambda)$ for each $i=1,2,3$.
 Let $\gamma\colon \Gamma\longrightarrow X_1\times_SX_2$
 (resp.\ $\gamma'\colon \Gamma'\longrightarrow X_2\times_SX_3$,
 resp.\ $\gamma''\colon \Gamma''\longrightarrow X_1\times_SX_3$)
 be an $S$-morphism such that $\gamma_1$ (resp.\ $\gamma'_1$, resp.\ $\gamma''_1$) is proper,  
 and $(\gamma_{\eta},c)$
 (resp.\ $(\gamma'_{\eta},c')$, resp.\ $(\gamma''_{\eta},c'')$)
 a cohomological correspondence from $L_1$ to $L_2$ (resp.\ from $L_2$ to $L_3$,
 resp.\ from $L_1$ to $L_3$).
 Moreover we assume that $\gamma_{1,s}^{-1}(Y_1)=\gamma_{2,s}^{-1}(Y_2)$, 
 $\gamma_{1,s}^{-1}(Z_1)=\gamma^{-1}_{2,s}(Z_2)$, 
 $\gamma'^{-1}_{1,s}(Y_2)=\gamma'^{-1}_{2,s}(Y_3)$, 
 $\gamma'^{-1}_{1,s}(Z_2)=\gamma'^{-1}_{2,s}(Z_3)$, 
 $\gamma''^{-1}_{1,s}(Y_1)=\gamma''^{-1}_{2,s}(Y_3)$
 and $\gamma''^{-1}_{1,s}(Z_1)=\gamma''^{-1}_{2,s}(Z_3)$. 
 Then, as above, the morphisms
 \begin{align*}
  R\Gamma_c(i^*j_*j^!R\psi(c))&\colon R\Gamma_c(X_{1,s},i_1^*R{j_1}_*Rj_1^!R\psi L_1)\longrightarrow R\Gamma_c(X_{2,s},i_2^*R{j_2}_*Rj_2^!R\psi L_2),\\
  R\Gamma_c(i^*j_*j^!R\psi(c'))&\colon R\Gamma_c(X_{2,s},i_2^*R{j_2}_*Rj_2^!R\psi L_2)\longrightarrow R\Gamma_c(X_{3,s},i_3^*R{j_3}_*Rj_3^!R\psi L_3),\\
  R\Gamma_c(i^*j_*j^!R\psi(c''))&\colon R\Gamma_c(X_{1,s},i_1^*R{j_1}_*Rj_1^!R\psi L_1)\longrightarrow R\Gamma_c(X_{3,s},i_3^*R{j_3}_*Rj_3^!R\psi L_3)
 \end{align*}
 are induced. Assume that the composite of $(\gamma_{\eta},c)$ and $(\gamma'_{\eta},c')$ coincides with
 $(\gamma_{\eta}'',c'')$. Then we have $R\Gamma_c(i^*j_*j^!R\psi(c'))\circ R\Gamma_c(i^*j_*j^!R\psi(c))=R\Gamma_c(i^*j_*j^!R\psi(c''))$.
\end{cor}

\begin{prf}
 By Proposition \ref{prop:specialization-indep}, we may replace $\gamma''$ by 
 $\Gamma\times_{X_2}\Gamma'\longrightarrow X_1\times_SX_3$. Then the equality is clear,
 since all the operations for cohomological correspondences are compatible with composition. 
\end{prf}

\def\cprime{$'$}
\providecommand{\bysame}{\leavevmode\hbox to3em{\hrulefill}\thinspace}
\providecommand{\MR}{\relax\ifhmode\unskip\space\fi MR }
\providecommand{\MRhref}[2]{%
  \href{http://www.ams.org/mathscinet-getitem?mr=#1}{#2}
}
\providecommand{\href}[2]{#2}

\bigbreak\bigbreak

\noindent Tetsushi Ito\par
\noindent Department of Mathematics, Faculty of Science, Kyoto University, Kyoto, 606--8502, Japan\par
\noindent E-mail address: \texttt{tetsushi@math.kyoto-u.ac.jp}

\bigbreak

\noindent Yoichi Mieda\par
\noindent Faculty of Mathematics, Kyushu University, 744 Motooka, Nishi-ku, Fukuoka, 819--0395, Japan\par
\noindent E-mail address: \texttt{mieda@math.kyushu-u.ac.jp}

\end{document}